\documentclass[11pt]{article}
\usepackage{graphicx}
\usepackage{amsmath,amssymb,amsthm,amsfonts}
\usepackage{amssymb}
\newtheorem{thm}{Theorem}[section]

\newtheorem{lem}{Lemma}[section]

\newtheorem{rem}{Remark}[section]
\theoremstyle{definition}

\usepackage{appendix}

\numberwithin{equation}{section}

\DeclareMathSymbol{\C}{\mathalpha}{AMSb}{"43}

\newcommand{\R}{{\mathbb{R}}}
\newcommand{\h}{{\mathcal{H}}}

\newcommand{\bsub}{\begin{subequations}}
\newcommand{\esub}{\end{subequations}$\!$}

\begin{document}

\title{A constrained minimization problem related to two coupled pseudo-relativistic Hartree equations
}
\author{Wenqing Wang
,\ Xiaoyu Zeng\ \ and Huan-Song Zhou\thanks{Corresponding author. \newline Email address: W.Q. Wang: wangwenqing\_1234@163.com; X.Y.Zeng: xyzeng@whut.edu.cn; H.S.Zhou: hszhou@whut.edu.cn. }\\
\ \\
\small
 Center for Mathematical Sciences and Department of Mathematics\\
 \small  Wuhan University of Technology, Wuhan 430070, China   }

\date{}

\smallbreak \maketitle

\begin{abstract}
We are concerned with the following constrained minimization problem:
$$e(a_{1},a_{2},\beta) := \inf\left\{E_{a_{1},a_{2},\beta}(u_{1},u_{2}): \|u_{1}\|_{L^{2}(\mathbb{R}^{3})} = \|u_{2}\|_{L^{2}(\mathbb{R}^{3})} = 1\right\},$$
where $E_{a_{1},a_{2},\beta}$ is the energy functional associated to two coupled pseudo-relativistic Hartree equations involving three parameters $a_{1}, a_{2}, \beta$ and two trapping potentials $V_1(x)$ and $V_2(x)$. In this paper, we obtain the existence of minimizers of $e(a_{1},a_{2},\beta)$ for possible $a_{1}, a_{2}$ and $\beta$ under suitable conditions on the potentials, which generalizes the results of the papers [16,17,18] in different senses.

\end{abstract}


\vskip 0.2truein

\section{Introduction}

Recently, Yang and Yang \cite{15} studied the following type of elliptic equation
\begin{equation}\label{1.1000}
\sqrt{-\Delta+m^{2}}\, u(x) + V(x)u(x) =\mu u(x) + a \,\big(|x|^{-1} \ast u^{2}(x)\big)u(x), \  \ x \in \mathbb{R}^{3},
\end{equation}
\noindent where $\mu$ and $a$ are real parameters, $m > 0,$ $V(x)$ is a trapping potential, $``\ast "$ denotes the convolution of two functions in $\mathbb{R}^{3},$ $\sqrt{-\Delta+m^{2}}$ is the so-called pseudo-relativistic operator, which is defined in the Sobolev space $
H^{\frac{1}{2}}(\mathbb{R}^{3})$ via Fourier transform, that is,
\begin{equation}\label{1.10000}
\sqrt{-\Delta+m^{2}}\, u(x) := \mathcal{F}^{-1}(\sqrt{|2\pi\xi|^{2} + m^{2}}\mathcal{F}u(\xi) )(x), \ \ u \in H^{\frac{1}{2}}(\mathbb{R}^{3}),
\end{equation}
where
$$\mathcal{F}u(\xi) = \int _{\R ^3} \exp^{-2\pi i \xi \cdot x} u(x)dx, \ \,  \ u \ \in \ H^{\frac{1}{2}}(\mathbb{R}^{3}) \text{ and }$$
$$H^{\frac{1}{2}}(\mathbb{R}^{3}) =  \{u\in L^{2}(\R ^3; \R): \|u\|^{2}_{H^{\frac{1}{2}}(\mathbb{R}^{3})} := \int _{\R ^3} (1 + 2\pi|\xi|)\Big|\mathcal{F}u(\xi)\Big|^{2} d\xi < \infty  \}.$$
In what follows, we denote by $\big(\cdot, \cdot\big)$ the complex $L^2$ inner product, i.e.,
$$ (u, v ) := \int _{\R ^3} u\bar{v} dx, ~~\forall~ u, v \in L^2 (\R ^3, \mathbb{C} ),$$
and then
\begin{equation}\label{1.1003400}
 (\sqrt{-\Delta+m^{2}}\, u, u ) := \int _{\R ^3}\sqrt{|2\pi\xi|^{2} + m^{2}}\Big|\mathcal{F}u(\xi)\Big|^{2}d\xi,~~~u \in H^{\frac{1}{2}}(\mathbb{R}^{3}).
\end{equation}
Moreover, for the norm of $H^{\frac{1}{2}}(\mathbb{R}^{3}),$ it is known that, see e.g., \cite{14, dd},
\begin{equation}\label{1.100000}
\|u\|^{2}_{H^{\frac{1}{2}}} = \int _{\R ^3}\left[u(-\Delta)^{\frac{1}{2}}u + u^{2}\right]dx = \int _{\R ^3}\left[|(-\Delta)^{\frac{1}{4}}u|^{2} + u^{2}\right]dx.
\end{equation}
Equation (\ref{1.1000}) occurs while investigating solitary waves (i.e., solutions in the form of $\psi(t, x) = e^{-i\mu t}u(x)$) for the following pseudo-relativistic Hartree equation:
$$i\partial_{t}\psi(t, x) = \sqrt{-\Delta+m^{2}}\psi(t, x) + V(x)\psi(t, x) - a(|x|^{-1}\ast|\psi|^{2})\psi(t, x), \
 x \  \in \ \mathbb{R}^{3}.$$
This equation is an important model for studying the dynamics of pseudo-relativistic boson stars, which has been studied particularly in the case of $V(x) \equiv -m,$ see e.g., \cite{0, 5, 3, 11, LY} and the references therein, for more backgrounds and related results. Motivated by the discussions on the Gross-Pitaevskii equation (see \cite{ee, 8, 10, 7}, etc.), the authors of paper \cite{15} studied equation (1.1) with $V(x) \not\equiv  constant$ and obtained $L^{2}$-normalized least energy solution by studying the constrained minimization problem:
\[
e(a):= \inf\left\{E_{a}(u): u \in \h ~\text{and} ~\|u\|_{2}^{2} := \int _{\R ^3} u^{2}(x)dx = 1\right\},
\]
where
$$\h :=  \{ u\in  H^{\frac{1}{2}}(\R ^3):\ \int _{\R ^3}  V(x)u^2 dx < \infty  \} \text{ and }$$
$$E_{a}(u) = \int_{\R ^3}  (u\sqrt{-\Delta+m^{2}}\, u + V(x)u^{2} )dx - \frac{a}{2}\int_{\mathbb{R}^{3}}(|x|^{-1}\ast u^{2})u^{2}dx.$$
Denote $a^{\ast} = \|Q\|_{2}^{2}$ with $Q > 0$ being a radially symmetric ground state solution (\cite{5, nn}) of
\begin{equation}\label{s0000}
\sqrt{-\Delta} u + u = \left(|x|^{-1}\ast u^2\right)u, \ x \ \in \ \mathbb{R}^{3},\ u \ \in \ H^{\frac{1}{2}}(\R ^3).
\end{equation}
We mention that the uniqueness on the ground state $Q$ of (\ref{s0000}) is still open, it is proved in \cite[Lemma 4.3]{15} that every ground state solution of (\ref{s0000}) has the same $L^{2}$-norm, that is, $a^{\ast} = \|Q\|_{2}^{2}$ is a fixed number for any ground state $Q$ of (\ref{s0000}). Based on this fact, in \cite{15}, the authors showed that $e(a)$ has a minimizer for any $a \in (0, a^{\ast}),$ but $e(a)$ has no minimizer if $a \geq a^{\ast}.$ This result implies that the single equation (1.1) has always a least energy solution with given $L^{2}$-norm as $a \in \big(0, a^\ast\big).$

In this paper, we are interested in the following system of two equations like (1.1) with coupled Hartree type nonlinearities, that is,
\begin{equation}\label{sys}
\begin{cases}
 \sqrt{-\Delta+m^{2}}\, u_{1}(x)+V_1(x)u_{1}=\mu_{1}u_{1}+a_1\left(|x|^{-1}\ast u_{1}^{2}\right)u_{1}
+\beta\left(|x|^{-1}\ast u_{2}^{2}\right)u_{1}, \\
 \sqrt{-\Delta+m^{2}}\, u_{2}(x)+V_2(x)u_{2}=\mu_{2}u_{2}+a_2\left(|x|^{-1}\ast u_{2}^{2}\right)u_{2}
+\beta\left(|x|^{-1}\ast u_{1}^{2}\right)u_{2}.
\end{cases}
\end{equation}
When (\ref{sys}) is a {\it Laplacian system}, that is, the pesudo-relativistic operator $\sqrt{-\Delta+m^{2}}$ is replaced by $-\Delta$ in (\ref{sys}), which was studied by Shi and Wang \cite{1} under $V_{1}(x) = V_{2}(x) \equiv 0,$ and they discussed how the existence of solutions in $H_{r}^{1}(\mathbb{R}^{3})\times H_{r}^{1}(\mathbb{R}^{3})$ depend on the parameters of $\mu_{1}, \mu_{2}, a_{1}, a_{2}$ and $\beta.$ Motivated by \cite{8}, where a Gross-Pitaevskii system (i.e., the Hartree nonlinear terms in the Laplacian system (\ref{sys}) are replaced by pure power type nonlinear terms) was investigated, Wang and Yang \cite{9} studied also the Laplacian system (\ref{sys}) with $V_{1}(x)$ and $V_{2}(x)$ being coercive functions, and they established some similar results to the paper \cite{8} under suitable ranges for the parameters $a_{1}, a_{2}$ and $\beta.$ In this paper, we aim to consider the system (\ref{sys}) with {\it both the pesudo-relativistic operator and the Hartree nonlinear terms} which are nonlocal, these make our discussions more complicated and difficult. Our results not only extend that of \cite{15} to the system case, but also extend that of \cite{1, 9} to the pesudo-relativistic operator system.
\vskip4pt
\indent In order to get solutions for the system (\ref{sys}), we may consider the following constrained variational minimization problem:
\begin{equation}\label{eq1.4}
\hat e(a_1,a_2,\beta):=\inf_{\{(u_1,u_2)\in \mathcal{M}\}} E_{a_1,a_2,\beta}(u_1,u_2),
\end{equation}
where
$\mathcal{M}=\Big\{(u, v)\in\mathcal{X}:\, \int _{\R ^3} u^2dx=\int _{\R ^3} v^2dx=1\Big\}$
and $\mathcal{X}=\h_1\times \h_2$ with
$$\h_i = \Big \{u \in  H^{\frac{1}{2}}(\R ^3):\ \int _{\R ^3}  V_i(x)u^2 dx < \infty \Big\} \  (\ i=1,\,2),  \text{ and }$$
\begin{equation}\label{eq1.789}
\|u\|^{2}_{_{\h_i}}= \int _{\R ^3} \left[|(-\Delta)^{\frac{1}{4}}u|^{2} + V_i(x)u^2\right]dx, \ \|u\|^{2}_{_\mathcal{X}}= \|u\|^{2}_{_{\h_1}}+\|u\|^{2}_{_{\h_2}}\ .
\end{equation}
Furthermore, it is not difficult to know that $\|\cdot\|_{_{\h_i}}$ is equivalent to the norm
$\|u\|_{m} =  (\int_{\R ^3}  [u\sqrt{-\Delta+m^{2}}\, u(x) + V(x)u^2 ]dx )^\frac{1}{2},$
 and the variational functional $E_{a_1,a_2,\beta}(u_1,u_2)$ defined  by
\begin{align}\label{f1}
   E_{a_1,a_2,\beta}(u_1,u_2)=&\sum_{i=1}^2\int_{\R ^3} \Big[u_i\sqrt{-\Delta+m^{2}}\, u_i(x) + V_i(x)u_i^2\Big]dx -\sum_{i=1}^2\frac{a_i}{2}\iint\limits_{\R ^3\times\R ^3} \frac{u_{i}^{2}(x)u_{i}^{2}(y)}{|x-y|}dxdy \nonumber\\
  &-\beta \iint_{\R ^3\times\R ^3} \frac{u_{1}^{2}(x)u_{2}^{2}(y)}{|x-y|}dxdy, \text{ for any }  (u_1,u_2) \in \mathcal{X}
\end{align}
is well defined and of $C^{1}$ in $\mathcal{X}$ by applying the following type of Gagliardo-Nirenberg inequality, for $ u \in H^{\frac{1}{2}}(\R ^3)$,
\begin{equation}\label{GNineq}
\int_{\R ^3}\Big(|x|^{-1}\ast |u|^{2}\Big)u^{2}dx =\iint_{\R ^3\times\R ^3}  \frac{u^{2}(x)u^{2}(y)}{|x-y|}dxdy \le \frac 2 {\|Q\|_2^{2}}  \|(-\Delta)^{\frac{1}{4}} u\|_{2}^2\|u\|^{2}_{2} .
\end{equation}
Particularly, ``=" holds in (\ref{GNineq}) if and only if $u(x) = Q(|x|)$ is a ground state of (\ref{s0000}), up to a translation and rescaling. Moreover, it follows from
(\ref{s0000}) and \eqref{GNineq} that
\begin{equation}\label{1:id}
 \|(-\Delta)^{\frac{1}{4}} Q \|_{2}^2  =\| Q \|_{2}^2=\frac{1}{2}\int_{\R ^3}\Big(|x|^{-1}\ast |Q|^{2}\Big)Q^{2}dx ,
\end{equation}
see \cite{5, nn} for more details.

\indent It is known that any minimizer $(u_{1}, u_{2}) \in \mathcal{M}$ of the minimization problem (\ref{eq1.4}) is a nontrivial solution (i.e., $u_{1}\not\equiv 0 \text{ and } u_{2}\not\equiv 0$)  for the system (\ref{sys}) with $\mu_{1}, \mu_{2}$ being suitable Lagrange multipliers. So, in this paper we focus on the study of the existence of minimizers for the minimization problem (\ref{eq1.4}). Motivated by  paper \cite{8} where a Laplacian system with power type nonlinearities was considered, we study problem (\ref{eq1.4}) by introducing a suitable auxiliary minimization problem, see (\ref{2:eq1.5}), upon which we can establish a criterion on the existence of minimizers for problem (\ref{eq1.4}), see our Theorem \ref{prop:2.1} in section 2. This criterion is crucial in proving our main results of the paper.
Furthermore,  in proving Theorem 1.3, we have to show that a rescaled minimizing  sequence of (\ref{eq1.4}) $\{ (\tilde{w}_{1n}, \tilde{w}_{2n}) \}$  has a nontrivial weak limit, that is, both components of the weak limit are nonzero,
and then the concentration-compactness principle are used. Since the sequence $\{ (\tilde{w}_{1n}, \tilde{w}_{2n}) \}$ has two components $\{\tilde{w}_{in}\} (i=1,2)$, if we apply the  concentration-compactness principle for $\{\tilde{w}_{in}\} $ separately,  we should get two sequence $\{y^{(i)}_{\epsilon_n}\} \subset \R^3  (i=1,2) $ accordingly. However, when we make translations for  $\{\tilde{w}_{in}\} (i=1,2)$ it is necessary to have  $\{y^{(i)}_{\epsilon_n}\}    (i=1,2) $  to be the same, which cannot be proved  and  the ideas of papers \cite{8,9} (for Laplacian system) do not work for our problem either. In this paper, by introducing some new strategies we avoid this kind of problem, see the proofs of  (\ref{5:scal})  for the details.

Before stating our main theorems of the paper, we need the following assumptions on the potentials $V_{i}(x)\ (i = 1, 2)$:
\vskip3pt
{\it
 {\bf ($V_{1}$)}~~
$V_i(x) \in L^\infty_{\rm
loc}(\R^3),\,\ \lim\limits_{|x|\to\infty} V_i(x) = \infty$ $\text{and}$ $\inf\limits_{x \in \R^3} V_i(x) = 0, \ i = 1, 2.$}

{\it
 {\bf ($V_{2}$)}~~
$\inf\limits_{x \in \R^3} \big (V_1(x) + V_2(x)\big)$ $\text{and}$ $\inf\limits_{x \in \R^{3}} V_i(x)(i = 1,2)$ are attained.}\\

\indent Now, we give the main theorems of the paper.
\begin{thm}\label{thm1}
Let $Q(x)$ be a positive radially symmetric ground state solution of
(\ref{s0000}) and $a^*=\|Q\|_2^2$. Under the condition $(V_{1}),$ we have
\begin{enumerate}
\item [\bf(i)] For any $a_{1}, a_{2} \in (0, a^{\ast})$,  if
\begin{equation}\label{gfd}
\beta  < \beta_{\ast} := \sqrt {(a^*-a_1)(a^*-a_2)},\end{equation}
then problem (\ref{eq1.4}) has a minimizer.

\item [ \bf(ii)] If $a_1 > a^*,$ or $a_2 > a^*,$ or $\beta > \beta^{\ast}$ with
\begin{equation}\label{ijbv}
\beta^{\ast}: = \frac{a^*-a_1}{2}+\frac{a^*-a_2}{2},\end{equation}
then problem (\ref{eq1.4}) has no minimizer.
\end{enumerate}
\end{thm}
\begin{rem}\label{rem1}
\begin{enumerate}
\item[\bf(i)] In Theorem 1.1 (i), $\beta \leq 0$ is allowed. In particular, if $\beta = 0,$ the system (\ref{sys}) is essentially equivalent to the single equation (\ref{1.1000}), and our existence result Theorem (i) is consistent with that of \cite{15}.
\item[ \bf(ii)] For $\beta_{\ast}$ and $\beta^{\ast}$ defined by (\ref{gfd}) and (\ref{ijbv}), it is easy to see that $\beta_{\ast} \leq \beta^{\ast}$ if
$a_{1}, a_{2} \in (0, a^{\ast}),$ but $\beta_{\ast} < \beta^{\ast}$ if $a_{1} \neq a_{2},$ and $\beta_{\ast} = \beta^{\ast}$ if $a_{1} = a_{2}.$
\end{enumerate}
\end{rem}

We mention that the above Theorem 1.1 does not work when
\begin{equation}\label{eq1.13}
a_{1}, a_{2} \in (0, a^{\ast}) \text{ and } \beta \in \big[\beta_{\ast}, \beta^{\ast}\big].
\end{equation}
In this case, we have the following theorem.
\begin{thm}\label{thm2} Under the conditions $(V_{1})$ and $(V_{2}),$ if $0 < a_1 \not= a_2 < a^*$ satisfies
  \begin{equation}\label{eq1.14}  |a_1 - a_2| \leq 2\sqrt{(a^*-a_1)(a^*-a_2)} = 2\beta_{\ast},\end{equation} then, there exists $\delta := \delta_{a_1, a_2} \in \big(\beta_{\ast}, \beta^{\ast}\big]$ such that problem (\ref{eq1.4}) has a minimizer for any $\beta \in \big[\beta_{\ast}, \delta\big).$
\end{thm}
It is easy to see that  Theorems 1.1 and 1.2 do not cover all the possibilities of the parameters $a_{1}, a_{2}$ and $\beta.$ To the best of authors' knowledge, the rest cases of Theorems 1.1 and 1.2 are still open. However, if $a_{1} = a_{2}$ in Theorem 1.2, we have only the following partial result.

\begin{thm}\label{thm03}
Under the conditions $(V_{1})$ and $(V_{2}),$ if $a_1 = a_2$ and (\ref{eq1.13}) is satisfied \big(i.e., $\beta =\beta_{\ast}=\beta^{\ast} =a^{\ast} - a_{1} = a^{\ast} - a_{2} > 0$\big), then
  \begin{enumerate}
    \item [\rm(i)] Problem (\ref{eq1.4}) has no minimizer provided $\underset{x \in \R^3}\inf\big(V_1(x) + V_2(x)\big) = 0$.
    \item [\rm(ii)]Problem (\ref{eq1.4}) has at least a minimizer provided
\begin{equation}\label{eq1.1400003}
\hat e(a^*-\beta,a^*-\beta,\beta) < \underset{x\in \R^3}\inf\big(V_1(x) + V_2(x)\big).
\end{equation}
\end{enumerate}
\end{thm}
We note that there are $V_{1}(x)$ and $V_{2}(x)$ such that condition (\ref{eq1.1400003}) holds.
For example, take $x_1, \ x_2 \in \R^3$ with $|x_1-x_2| > 5$ and $\zeta_i(x) \in C^2_0(B_1(x_i), \R^+)$ with $\|\zeta_i\|_2 = 1$,
let $V_i(x)\in C^2(\R^3, \R^{+})$ with $\lim\limits _{|x|\to\infty}V_i(x)=\infty\  (i = 1,\,2)$ such that
$$V_i(x)|_{ B_1(x_i)} =0 \text{ and } V_i (x)_{B_2^c(x_i)}\ge 2c_\zeta, $$
\text{ where }  $c_\zeta \in (0, +\infty)$ \text{ is given by}
$$
c_\zeta=\sum_{i=1}^2\int_{\R ^3} \Big(\zeta_i\sqrt{-\Delta + m^2} \, \zeta_i-\frac{a^*-\beta}{2}\phi_{\zeta_i}\zeta_i^2\Big)dx- \beta\int_{\R ^3}\phi_{\zeta_1}\zeta_2^2dx \text{ with } \phi_{\zeta_1}(x) =|x|^{-1}\ast \zeta_1(x)^{2}.
$$
By these facts, we see that
$$ \hat e(a^*-\beta,a^*-\beta,\beta)\leq E_{a^*-\beta,a^*-\beta,\beta}(\zeta_1,\zeta_2)=c_\zeta<2c_\zeta\le \inf_{x\in\R^3} \big(V_1(x)+V_2(x)\big).$$

The paper is organized as follows. In Section \ref{Sec2}, an auxiliary minimization problem is introduced, and some properties about it are proved, which are useful in showing our main Theorems. In Section \ref{Sec3}, we give the proofs of Theorems \ref{thm1}-\ref{thm03}.

\section{An auxiliary minimization problem and related results}\label{Sec2}

In this section, we focus on establishing a criterion on the existence of minimizers for problem (\ref{eq1.4}) by introducing an auxiliary constrained minimization problem. For this purpose, we first recall some useful lemmas.

\begin{lem}(\cite[Lemma 2.1]{15})  \label{2:lem1}
Suppose $V_{i} \in L_{loc}^\infty(\R^3)$ and
$\lim\limits_{|x|\to\infty}V_i(x)=\infty, i = 1, 2.$ Then the embedding $\h_1, \h_2\hookrightarrow L^{q}(\R^3)$ is compact for $2 \leq q < 3.$
\end{lem}
\vskip2pt

\begin{lem}\label{lemA2}(\cite[Theorem 4.3]{14})
Let $f \in L^{p}(\mathbb{R}^3)$ and $g \in L^{q}(\mathbb{R}^3).$ Then,
$$\iint_{\R ^3\times\R ^3} \frac{f(x)g(y)}{|x-y|^t}dxdy \leq c(p,q,t)\|f\|_{p}\|g\|_{q},$$
where $1< p, q < \infty, 0 < t < 3$ and $\frac{1}{p} + \frac{1}{q} + \frac{t}{3} =2.$  ~\hfill$\square$
\end{lem}

\begin{lem}\label{lemA3}(\cite[Lemma 2.2]{1})
Let $u, v \in L^{\frac{12}{5}}(\mathbb{R}^3)$. Then,
$$\iint_{\R ^3\times\R ^3}\frac{u^{2}(x)v^{2}(y)}{|x-y|}dxdy \leq \left(\iint_{\R ^3\times\R ^3}\frac{u^{2}(x)u^{2}(y)}{|x-y|}dxdy\iint_{\R ^3\times\R ^3} \frac{v^{2}(x)v^{2}(y)}{|x-y|}dxdy\right)^{\frac{1}{2}}
,$$
and the equality holds if and only if $u^{2}(x) = \kappa v^{2}(x)$ for some $\kappa > 0.$ ~\hfill$\square$
\end{lem}
For $a_{1}, a_{2} > 0$ and $\beta \in \R,$ we introduce the following minimization problem
\begin{equation}\label{2:eq1.5}
\eta(a_1,a_2,\beta):=\inf\big\{J_{a_{1},a_{2},\beta} (u_{1},u_{2} ):  u_{i}  \in   H^{\frac{1}{2}}(\R ^3), \|u_i\|_2^2=1, ~i=1,2\big\},
\end{equation}
where
\begin{equation}\label{5}
J_{a_{1},a_{2},\beta}\big(u_{1},u_{2}\big) = \frac{2\|(-\Delta)^{\frac{1}{4}} u_{1}(x) \|_{2}^2+2\|(-\Delta)^{\frac{1}{4}} u_{2}(x) \|_{2}^2}{a_1\int_{\R ^3}\phi_{u_1}u_{1}^{2}dx+a_2\int_{\R ^3}\phi_{u_2}u_{2}^{2}dx+2\beta^{+}\int_{\R ^3}\phi_{u_1}u_{2}^{2}dx}, \text{ and }
\end{equation}
$$\phi_{u}(x) :=|x|^{-1} \ast u^2(x)= \int_{\R ^3}\frac{u^{2}(y)}{|x-y|}dy, ~~\beta^{+} = \max\big\{0, \beta\big\}.$$
Then, we have
\begin{lem}\label{lemA4}
$\eta(a_1,a_2,\beta)$ defined by (\ref{2:eq1.5}) is locally Lipschitz continuous with respect to $(a_1,a_2,\beta) \in \R^3_+ \backslash ~\big\{(0,0,0)\big\}$.
\end{lem}

\noindent\textbf{Proof.}  Motivated by \cite[Lemma 2.2]{8},  for $Q(x)$ being given by \eqref{1:id}, we take $(u_{1}, u_{2}) = \Big(\frac{Q(x)}{\|Q\|_2},\frac{Q(x)}{\|Q\|_2}\Big)$ in (\ref{5}), then it follows from (\ref{1:id}) and (\ref{2:eq1.5}) that
\begin{equation}\label{sqdsz}\eta(a_1,a_2,\beta) \leq \frac{2a^*}{a_1+a_2+2\beta}.\end{equation}
Moreover, by Gagliardo-Nirenberg inequality (\ref{GNineq}) and Lemma 2.3, for any $(u_1, u_2)\in H^\frac{1}{2}(\R^3)\times H^\frac{1}{2}(\R^3)$ satisfying $\int_{\mathbb{R}^{3}} |u_i|^2dx=1, i=1,\,2,$ we obtain that, for any $(a_1,a_2,\beta) \in \R^3_+ \backslash ~\big\{(0,0,0)\big\},$
\begin{align*}
J_{a_{1},a_{2},\beta}\big(u_{1},u_{2}\big) &= \frac{2\|(-\Delta)^{\frac{1}{4}} u_{1}(x) \|_{2}^2+2\|(-\Delta)^{\frac{1}{4}} u_{2}(x) \|_{2}^2}{a_1\int_{\R ^3}\phi_{u_1}u_{1}^{2}dx+a_2\int_{\R ^3}\phi_{u_2}u_{2}^{2}dx+2\beta\int_{\R ^3}\phi_{u_1}u_{2}^{2}dx}\\
&\geq \frac{a^*\Big(\int_{\R ^3}\phi_{u_1}u_{1}^{2}dx + \int_{\R ^3}\phi_{u_2}u_{2}^{2}dx\Big)}{(a_1+\beta)
\int_{\R ^3}\phi_{u_1}u_{1}^{2}dx+(a_2+\beta)\int_{\R ^3}\phi_{u_2}u_{2}^{2}dx}\geq \frac{a^*}{\max\{a_1+\beta,a_2+\beta\}}.
\end{align*}
Hence,
\begin{equation}\label{eq7.2}
\frac{a^*}{\max\{a_1+\beta,a_2+\beta\}}\leq \eta(a_1,a_2,\beta) \leq \frac{2a^*}{a_1+a_2+2\beta}.
\end{equation}

 Let $\{(u_{1n},u_{2n})\} \subset H^\frac{1}{2}(\R^3)\times H^\frac{1}{2}(\R^3)$ be a minimizing sequence of $\eta(a_1,a_2,\beta),$ that is,
\begin{equation}\label{equation56}
\lim_{n \rightarrow \infty} J_{a_{1},a_{2},\beta}\big(u_{1n},u_{2n}\big) = \eta(a_1,a_2,\beta),~ \|u_i\|_2^2=1,i=1,2.\end{equation}
Since (\ref{5}) is invariant under the rescaling: $u(x)\mapsto \lambda^{\frac{3}{2}} u(\lambda x), \lambda>0$,  we may assume that
\begin{equation}\label{eq7.3}
 \|(-\Delta)^{\frac{1}{4}} u_{1n}(x) \|_{2}^2+\|(-\Delta)^{\frac{1}{4}} u_{2n}(x) \|_{2}^2=1\  \text{ for all } n\in\mathbb{N}^+.
\end{equation}
Then, the Gagliardo-Nirenberg inequality (\ref{GNineq}) and Lemma 2.3 imply that
$$\int_{\R ^3}\phi_{u_{in}}u_{in}^{2}dx\leq \frac{2}{a^*}\|(-\Delta)^{\frac{1}{4}} u_{in}(x) \|_{2}^2\le \frac{2}{a^*},\, \ i=1,\,2, $$
and
\begin{align*}\int_{\R ^3}\phi_{u_{1n}}u_{2n}^{2}dx&\leq \frac{1}{2}\Big(\int_{\R ^3}\phi_{u_{1n}}u_{1n}^{2}dx + \int_{\R ^3}\phi_{u_{2n}}u_{2n}^{2}dx\Big)\\
&\leq
\frac{1}{a^*} (\|(-\Delta)^{\frac{1}{4}} u_{1n}(x) \|_{2}^2+\|(-\Delta)^{\frac{1}{4}} u_{2n}(x) \|_{2}^2 )=\frac{1}{a^*}.\end{align*}
In addition, for any $(a_1,a_2,\beta),(\tilde a_1, \tilde a_2,\tilde\beta)\in \R^3_+ \backslash ~\big\{(0,0,0)\big\},$ it follows from (\ref{equation56}) and (\ref{eq7.3}) that
\begin{align*}
\frac{1}{\eta(a_1,a_2,\beta)}
=&\lim_{n\to\infty}\Big[\frac{\tilde{a}_1\int_{\R ^3}\phi_{u_{1n}}u_{1n}^{2}dx+\tilde{a}_2\int_{\R ^3}\phi_{u_{2n}}u_{2n}^{2}dx+2\tilde{\beta}\int_{\R ^3}\phi_{u_{1n}}u_{2n}^{2}dx}{2\|(-\Delta)^{\frac{1}{4}} u_{1n}(x) \|_{2}^2+2\|(-\Delta)^{\frac{1}{4}} u_{2n}(x) \|_{2}^2}\\
 & \hspace{1.5cm}+\overset{2}{\sum_{i=1}}\frac{a_i-\tilde a_i}{2}\int_{\R ^3}\phi_{u_{in}}u_{in}^{2}dx+(\beta-\tilde \beta)\int_{\R ^3}\phi_{u_{1n}}u_{2n}^{2}dx \Big]\\
\leq &\lim_{n\to\infty}\Big[\frac{1}{\eta(\tilde a_1,\tilde a_2,\tilde\beta)}
+\overset{2}{\sum_{i=1}}\frac{|a_i-\tilde a_i|}{2}\int_{\R ^3}\phi_{u_{in}}u_{in}^{2}dx+|\beta-\tilde \beta|\int_{\R ^3}\phi_{u_{1n}}u_{2n}^{2}dx\Big] \\
\leq &\frac{1}{\eta(\tilde a_1,\tilde a_2,\tilde\beta)}
+\overset{2}{\sum_{i=1}}\frac{|a_i-\tilde a_i|}{a^*}+\frac{|\beta-\tilde \beta|}{a^*},
\end{align*}
i.e.,
\begin{equation}\label{eq7.4}
\frac{1}{\eta(a_1,a_2,\beta)}-\frac{1}{\eta(\tilde a_1,\tilde a_2,\tilde\beta)}\leq \frac{3}{a^*}\big|(a_1,a_2,\beta)-(\tilde a_1, \tilde a_2,\tilde\beta)\big|.
\end{equation}
Similarly, for a minimizing sequence $\{(\tilde u_{1n},\tilde u_{2n})\}$ of
$\eta(\tilde a_1,\tilde a_2,\tilde\beta)$, we have
\begin{equation}\label{saft}
\frac{1}{\eta(\tilde a_1,\tilde a_2,\tilde\beta)}-\frac{1}{\eta(a_1,a_2,\beta)}\leq \frac{3}{a^*}\big|(a_1,a_2,\beta)-(\tilde a_1, \tilde a_2,\tilde\beta)\big|.\end{equation}
Hence, (\ref{eq7.4}) and (\ref{saft}) imply that
\begin{align*}\Big|\frac{\eta(a_1,a_2,\beta)-\eta(\tilde a_1,\tilde a_2,\tilde\beta)}{\eta(a_1,a_2,\beta)\eta(\tilde a_1,\tilde a_2,\tilde\beta)}\Big|
&=\Big|\frac{1}{\eta(a_1,a_2,\beta)}-\frac{1}{\eta(\tilde a_1,\tilde a_2,\tilde\beta)}\Big|
\leq \frac{3}{a^*}\big|(a_1,a_2,\beta)-(\tilde a_1, \tilde a_2,\tilde\beta)\big|,\end{align*}
this and (\ref{eq7.2}) give that
$$\big|\eta(a_1,a_2,\beta)-\eta(\tilde a_1,\tilde a_2,\tilde\beta)\big|\leq \frac{12a^*}{(a_1+a_2+2\beta)(\tilde a_1+\tilde a_2+2\tilde\beta)}\big|(a_1,a_2,\beta)-(\tilde a_1, \tilde a_2,\tilde\beta)\big|.$$
So, $\eta(\cdot)$ is locally Lipschitz continuous in $\R^3_+ \backslash ~\big\{(0,0,0)\big\}$.  \hfill$\square$

\begin{rem} Since $\eta(a_1,a_2,\beta) \equiv \eta(a_1,a_2,0)$ if $\beta \leq 0,$ it is easy to see that Lemma \ref{lemA4} is also true for $(a_{1},a_{2}) \in \R^2_+ \backslash ~\big\{(0,0)\big\}$ and $\beta \leq 0.$ \end{rem}

\begin{lem}\label{33:lem3.1}
Let $a_1, a_2 \in (0, a^*)$ and $a_1 \neq a_2.$ Then\\
\indent {\bf(i)} $\eta\big(a_1, a_2, \beta_{\ast}\big) > 1$ if $|a_1 - a_2|\leq 2\beta_{\ast} = 2\sqrt {(a^*-a_1)(a^*-a_2)};$\\
\indent {\bf(ii)} $\eta\big(a_1, a_2, \beta^{\ast}\big) \leq 1,$ $\beta^{\ast}$ is given by (\ref{ijbv}).
\end{lem}

\noindent\textbf{Proof.}
$(\bf i)$ Let $\{(u_{1n},u_{2n})\}$ be a  minimizing sequence of $\eta\big(a_1, a_2, \beta_{\ast}\big).$ Then, we may assume that $\{(u_{1n},u_{2n})\}\subset H_r^\frac{1}{2}(\R^3)\times H_r^\frac{1}{2}(\R^3)$ (i.e., $H^\frac{1}{2}(\R^3)$ with radially symmetric functions) due to the property of Schwartz symmetrization, see e.g., \cite[Theorem 3.7]{14}.
  Similar to (\ref{eq7.3}), we suppose that
\begin{equation}\label{eq7.305}
 \|(-\Delta)^{\frac{1}{4}} u_{1n}(x) \|_{2}^2+\|(-\Delta)^{\frac{1}{4}} u_{2n}(x) \|_{2}^2 = 1\  \text{ for all } n\in\mathbb{N}^+.
\end{equation}
This and (\ref{eq7.2}) as well as the Gagliard-Nirenberg inequality (\ref{GNineq}) show that, passing to a subsequence if necessary,
\begin{equation}\label{eqeqeq}
  0 < C_{0} \leq a_1\int_{\R ^3}\phi_{u_{1n}}u_{1n}^{2}dx+a_2\int_{\R ^3}\phi_{u_{2n}}u_{2n}^{2}dx+2\beta\int_{\R ^3}\phi_{u_{1n}}u_{2n}^{2}dx \leq \bar{C}_{0} < \infty.
\end{equation}
Now, we prove $(\bf i)$ by considering the following {\bf two cases}:
 \vskip 0.1truein
 \noindent \textbf{Case 1.} $\displaystyle \int_{\mathbb{R}^{3}} \phi_{u_{1n}}u_{1n}^{2}dx \rightarrow 0,$ or $\displaystyle \int_{\mathbb{R}^{3}} \phi_{u_{2n}}u_{2n}^{2}dx \rightarrow 0 ~\text{as}~ n \to \infty.$

  In this case, we assume that $\int_{\mathbb{R}^{3}} \phi_{u_{1n}}u_{1n}^{2}dx \rightarrow 0 ~\text{as}~ n \to \infty$  (the
other case can be done similarly). By $\int_{\mathbb{R}^{3}} \phi_{u_{1n}}u_{1n}^{2}dx \overset{n}\to 0$ and Lemma \ref{lemA3}, we know from \eqref{eqeqeq} that
$$\int_{\mathbb{R}^{3}} \phi_{u_{1n}}u_{2n}^{2}dx \overset{n}\to 0~ \text{and}~ \int_{\mathbb{R}^{3}} \phi_{u_{2n}}u_{2n}^{2}dx \geq C > 0.$$
 Then, by \eqref{GNineq} we have
 \begin{align}\eta(a_1,a_2,\beta_{\ast})
&=\lim_{n\to\infty} \frac{2\|(-\Delta)^{\frac{1}{4}} u_{1n}(x) \|_{2}^2+2\|(-\Delta)^{\frac{1}{4}} u_{2n}(x) \|_{2}^2}{a_2\int_{\mathbb{R}^{3}} \phi_{u_{2n}}u_{2n}^{2}dx+o(1)}\nonumber\\
&\geq \lim_{n\to\infty} \frac{a^{\ast}\int_{\mathbb{R}^{3}}\phi_{u_{2n}}u_{2n}^{2}dx}{a_2\int_{\mathbb{R}^{3}}\phi_{u_{2n}}u_{2n}^{2}dx+o(1)}\geq \frac{a^*}{a_2} > 1, ~\text{since}~ a_{2}  \in (0, a_{\ast}).\label{eq3.18}
\end{align}
Hence, $(\bf i)$ is proved.

 \vskip 0.1truein
\noindent \textbf{Case 2.} $\displaystyle\int_{\mathbb{R}^{3}} \phi_{u_{1n}}u_{1n}^{2}dx \geq C > 0$ and $\displaystyle \int_{\mathbb{R}^{3}} \phi_{u_{2n}}u_{2n}^{2}dx \geq C > 0$ for some $C>0$.

  In this case, it follows  from (\ref{GNineq}), (\ref{eq7.305}) and (\ref{eqeqeq}) that, there is a constant $C_1>0$ such that
\begin{equation}\label{eq7.5}
 C_{1} \leq \|(-\Delta)^{\frac{1}{4}} u_{1n}(x) \|_{2}^2, \|(-\Delta)^{\frac{1}{4}} u_{2n}(x) \|_{2}^2, \int_{\mathbb{R}^{3}} \phi_{u_{1n}}u_{1n}^{2}dx, \int_{\mathbb{R}^{3}} \phi_{u_{2n}}u_{2n}^{2}dx \leq \bar{C}_{1}.
\end{equation}
\noindent Since $\{u_{in}\} \subset H_r^\frac{1}{2}(\R^3)$ ($i=1,\,2$), by (\ref{eq7.5}) and the compactness of Sobolev embedding, we deduce that there exists $u_i(x)\in H_r^\frac{1}{2}(\R^3)$ such that
\begin{equation}\label{eq3.02}
u_{in} \overset{n}\rightharpoonup u_i  \text{ weakly in }    H_{r}^\frac{1}{2}(\R^3), \
u_{in} \overset{n}\to u_i    \text{ strongly in }   L^p(\R^3),   \   p\in(2,3),\  i=1, 2.
\end{equation}
By applying  Lemma 2.2 of \cite{15}, we know that
$$\int_{\R ^3}\phi_{u_in}u_{in}^{2}dx \overset{n}\to \int_{\R ^3}\phi_{u_i}u_{i}^{2}dx, ~~i~ = ~1,2.$$
Then,  (\ref{eq7.5}) implies that
\begin{equation}\label{agyt}\int_{\mathbb{R}^{3}} \phi_{u_i}u_{i}^{2}dx \geq C > 0 \text{ and }u_i\not\equiv 0, \ i=1,\,2.\end{equation}

Since $0 < a_1\not= a_2 < a^*$, without loss of generality, we may suppose that $a_1 < a_2$, and using our condition $|a_{1} - a_{2}| \leq 2\beta_{\ast} = 2\sqrt{(a^{\ast} - a_{1})(a^{\ast} - a_{2})},$ we have
\begin{equation}\label{eq3.03}
0 < a_1 < a_2 <a^* \ \text{ and } \ a_2 \leq 2\beta_{\ast} + a_1, ~\beta_{\ast} = \sqrt{(a^*-a_1)(a^*-a_2)} > 0.
\end{equation}
So, it follows from  Gagliardo-Nirenberg inequality (\ref{GNineq}), Lemma 2.3 and (\ref{eq3.02}) - (\ref{agyt})  that

\begin{align}
\eta(a_1,a_2,\beta_{\ast})
& =\lim_{n \rightarrow \infty} J_{a_{1},a_{2}, \beta_{\ast}} \big(u_{1n}, u_{2n} \big)  \nonumber  \\
& \geq \lim_{n\to\infty}\frac{a^* (\int_{\R ^3}\phi_{u_{1n}}u_{1n}^{2}dx + \int_{\R ^3}\phi_{u_{2n}}u_{2n}^{2}dx )}{a_1\int_{\R ^3}\phi_{u_{1n}}u_{1n}^{2}dx+a_2\int_{\R ^3}\phi_{u_{2n}}u_{2n}^{2}dx+2\beta_{\ast}\int_{\R ^3}\phi_{u_{1n}}u_{2n}^{2}dx} \nonumber \\
& =\frac{a^* (\int_{\R ^3}\phi_{u_1}u_{1}^{2}dx+\int_{\R ^3}\phi_{u_2}u_{2}^{2}dx )}{a_1\int_{\R ^3}\phi_{u_1}u_{1}^{2}dx+a_2\int_{\R ^3}\phi_{u_2}u_{2}^{2}dx+2\beta_{\ast}\int_{\R ^3}\phi_{u_1}u_{2}^{2}dx}\nonumber\\
& \geq \frac{a^* (\int_{\R ^3}\phi_{u_1}u_{1}^{2}dx+\int_{\R ^3}\phi_{u_2}u_{2}^{2}dx )}{a_1\int_{\R ^3}\phi_{u_1}u_{1}^{2}dx+a_2\int_{\R ^3}\phi_{u_2}u_{2}^{2}dx+2\beta_{\ast} (\int_{\R ^3}\phi_{u_1}u_{1}^{2}dx\int_{\R ^3}\phi_{u_2}u_{2}^{2}dx )^{\frac{1}{2}}}\label{2.17}\\
& =\gamma_{a_1,a_2,\beta_{\ast}}(s_{0})~~ \text{with}~~ s_{0}:=\Big(\frac{\int_{\R ^3}\phi_{u_2}u_{2}^{2}dx}{\int_{\R ^3}\phi_{u_1}u_{1}^{2}dx}\Big)^{\frac{1}{2}} > 0\nonumber\\
& \geq\inf_{t\in(0,\infty)}\gamma_{a_1,a_2,\beta_{\ast}}(t), \label{eq3.4}
\end{align}
where
\begin{equation}\label{eq2.6}
\gamma_{a_1,a_2,\beta_{\ast}}(t):= \frac{a^*(1+t^2)}{a_1+a_2t^2+2\beta_{\ast} t}.
\end{equation}
By Lemma 2.3, the inequality  (\ref{2.17}) becomes equality if and only if
\begin{equation}\label{eq3.6}
u_2^2(x)=\kappa u_1^2(x) \,\ \text{ for some }\, \ \kappa>0.
\end{equation}
For $\beta_{\ast} = \sqrt{(a^*-a_1)(a^*-a_2)},$ it is easy to see from the definition (\ref{eq2.6}) that
\begin{equation}\label{eq3.07}
\gamma_{a_1,a_2,\beta_{\ast}}(t) \geq 1,\,\  \forall\  t \in (0, +\infty),
\end{equation}
and
\begin{equation}\label{eq3.8}
\gamma_{a_1,a_2,\beta_{\ast}}(t) = 1 \,\ \Leftrightarrow\,\  t=t_0:=\sqrt\frac{a^*-a_1}{a^*-a_2}.
\end{equation}
Then, $\eta\big(a_1, a_2, \beta_{\ast}\big)\geq 1$ follows from (\ref{eq3.4}) and (\ref{eq3.07}). Therefore, to prove $(\bf i)$ we need only to rule out $\eta\big(a_1, a_2, \beta_{\ast}\big) = 1.$
In fact, if
\begin{equation}\label{xcvvv}
\eta\big(a_1, a_2, \beta_{\ast}\big) = 1, \end{equation}
then (\ref{eq3.6}) holds, that is
\begin{equation}\label{eq3.9}
u_2^2(x)=\kappa u_1^2(x),\,\ \text{and}\,\ \kappa = t_0 = \sqrt\frac{a^*-a_1}{a^*-a_2} > 1.
\end{equation}
Hence, it follows from (\ref{eq3.02}) and (\ref{eq3.9}) that
\begin{align}
\eta(a_1,a_2,\beta_{\ast})=&\lim_{n \rightarrow \infty} J_{a_{1},a_{2},\beta_{\ast}}\big(u_{1n},u_{2n}\big)\nonumber\\
&\geq \frac{2\|(-\Delta)^{\frac{1}{4}} u_{1}(x) \|_{2}^2+2\|(-\Delta)^{\frac{1}{4}} u_{2}(x) \|_{2}^2}{a_1\int_{\R ^3}\phi_{u_1}u_{1}^{2}dx+a_2\int_{\R ^3}\phi_{u_2}u_{2}^{2}dx+2\beta_{\ast}\int_{\R ^3}\phi_{u_1}u_{2}^{2}dx} \nonumber\\
&=\frac{2(1+\kappa)}{a_1+a_2\kappa^2+2\beta_{\ast}\kappa}\cdot\frac{\|(-\Delta)^{\frac{1}{4}} u_{1}(x) \|_{2}^2}{\int_{\R ^3}\phi_{u_1}u_{1}^{2}dx}.\label{eq3.10}
\end{align}

On the other hand, let
\begin{equation}\label{eq3.11}
\tilde u_i(x)=\frac{1}{\sqrt\lambda_i}u_i(x)\,\ \text{with}\,\  \lambda_i:=\int_{\R ^3} |u_i|^2dx\leq \lim_{n\to\infty}\int_{\R ^3} |u_{in}|^2dx=1.
\end{equation}
Then, $\int_{\R ^3}|\tilde u_i|^2dx = 1, ~~i=1,2,$ and
\begin{equation}\label{eq3.12}
\lambda_2 = \kappa\lambda_1 \leq 1, ~\text{by}~ (\ref{eq3.9}).
\end{equation}
Therefore, we know from (\ref{eq3.9}), (\ref{eq3.11}) and (\ref{eq3.12}) that
\begin{align}
\eta(a_1,a_2,\beta_{\ast})&\leq \frac{2\|(-\Delta)^{\frac{1}{4}} \tilde{u}_{1}(x) \|_{2}^2+2\|(-\Delta)^{\frac{1}{4}} \tilde{u}_{2}(x) \|_{2}^2}{a_1\int_{\R ^3}\phi_{\tilde{u}_1}\tilde{u}_{1}^{2}dx+a_2\int_{\R ^3}\phi_{\tilde{u}_2}\tilde{u}_{2}^{2}dx+2\beta_{\ast}\int_{\R ^3}\phi_{\tilde{u}_1}\tilde{u}_{2}^{2}dx}\nonumber\\
&=\frac{4\lambda_1}{a_1+a_2+2\beta_{\ast}}\cdot \frac{\|(-\Delta)^{\frac{1}{4}} u_{1}(x) \|_{2}^2}{\int_{\R ^3}\phi_{u_1}u_{1}^{2}dx}\nonumber
\leq \frac{4/\kappa}{a_1+a_2+2\beta_{\ast}}\cdot\frac{\|(-\Delta)^{\frac{1}{4}} u_{1}(x) \|_{2}^2}{\int_{\R ^3}\phi_{u_1}u_{1}^{2}dx}.\nonumber
\end{align}
Combing (\ref{eq3.10}), we see that
$$\frac{2(1+\kappa)}{a_1+a_2\kappa^2+2\beta_{\ast}\kappa}\leq \frac{4/\kappa}{a_1+a_2+2\beta_{\ast}},$$
i.e.,
\begin{equation}\label{xcv}
\frac{2(\kappa+\kappa^2)}{a_1+a_2\kappa^2+2\beta_{\ast}\kappa} \leq \frac{4}{a_1+a_2+2\beta_{\ast}}, \ \text{where }\ \kappa = \sqrt\frac{a^*-a_1}{a^*-a_2} > 1.\end{equation}
Now, we claim that (\ref{xcv}) can not be true, which then implies that (\ref{xcvvv}) is false, and $(\bf i)$ is proved.\\
\indent In fact, by (\ref{eq3.03}), we know that $\big(a^{\ast} - a_{1}\big) - \big(a^{\ast} - a_{2}\big) \leq 2\sqrt{(a^*-a_1)(a^*-a_2)},$
that is,
$$\frac{a^{\ast} - a_{1}}{a^{\ast} - a_{2}} - 1 \leq 2\sqrt\frac{a^*-a_1}{a^*-a_2},~ \Longleftrightarrow~ \kappa^{2} - 1 \leq 2\kappa,$$
hence
$1<\kappa \leq 1 + \sqrt{2}.$
Let
$$g(t) := \frac{2(t + t^2)}{a_1+a_2t^2 + 2\beta_{\ast} t}.$$
Then, $g(t)$ is a strictly increasing function as $t \in \big[1, 1 + \sqrt{2}\big],$ because
$$g^{'}(t) = 2\frac{(2\beta_{\ast} - a_{2})t^{2} + 2a_{1}t + a_{1}}{\big(a_{1} + a_{2}t^{2} + 2\beta_{\ast} t\big)^{2}},$$
and by (\ref{eq3.03}),
 $$(2\beta_{\ast} - a_{2})t^{2} + 2a_{1}t + a_{1} \geq -a_{1}t^{2} + 2a_{1}t + a_{1} = a_{1}\big[2 - (t - 1)^{2}\big],$$
which means that $g'(t) > 0$ for $ t \in \big[1, 1 + \sqrt{2}\big).$ So, $g(1) < g(\kappa)$ by $\kappa \in \big(1, 1 + \sqrt{2}\big],$ and then (\ref{xcv}) can not be true.

($\bf ii$): By (\ref{sqdsz}), for any $(a_1,a_2,\beta) \in \R^3_+ \backslash ~\big\{(0,0,0)\big\},$ we have
$$\eta(a_1,a_2,\beta) \leq \frac{2a^*}{a_1+a_2+2\beta}.$$
Take   $\beta = \beta^{\ast} = \frac{2a^{\ast}-a_{1}-a_{2}}{2} $ in the above inequality,  then it is easy to see that $\beta^{\ast}>0$ and $\eta(a_1,a_2,\beta^{\ast}) \leq 1$ since $a_1, a_2 \in (0, a^\ast)$.~\hfill$\square$

Finally we give our main result of this section.
\begin{thm}\label{prop:2.1} If $a_1, a_2 > 0,$ $\beta \in \mathbb{R}$ and the condition $(V_{1})$ holds. Then
\begin{enumerate}
\item [\bf (i)]  problem (\ref{eq1.4}) has a minimizer if $\eta(a_1,a_2,\beta) > 1.$
\item [\bf(ii)]  problem  (\ref{eq1.4}) has no minimizer if $\eta(a_1,a_2,\beta) < 1.$
\end{enumerate}
\end{thm}

\noindent\textbf{Proof.} ($\bf{i}$) Let
$\{(u_{1n},u_{2n})\}\subset \mathcal{M}$ be a minimizing sequence
of problem (\ref{eq1.4}),  then
$$\|u_{1n}\|_2^2=\|u_{2n}\|_2^2=1\quad \text{and}\quad
\lim_{n\to\infty}E_{a_1,a_2,\beta}(u_{1n},u_{2n})=\hat e(a_1,a_2,\beta).$$
Since
$\sqrt{-\Delta + m^{2}} \geq \sqrt{-\Delta}$,  i.e.,
$$ (\sqrt{-\Delta + m^{2}}u, u ) \geq  (\sqrt{-\Delta}u, u ) \text{ for } u \in H^{\frac{1}{2}}(\mathbb{R}^{3}),$$
we know from the definition of $\eta(a_1,a_2,\beta)$ in (\ref{2:eq1.5})   that
\begin{eqnarray}\label{eq2.1}
\begin{split}
  E_{a_1,a_2,\beta}(u_{1n},u_{2n})  \geq &\sum_{i=1}^2\|(-\Delta)^{\frac{1}{4}}u_{in}\|_{2}^{2}+ \int_{\R ^3} V_{i}(x)|u_{in}|^2dx\\
   & -\frac{a_i}{2}\int_{\R ^3} (|x|^{-1}\ast |u_{in}|^{2} )u_{in}^{2}dx  -\beta\iint_{\R ^3\times\R ^3}  \frac{u_{1n}^{2}(x)u_{2n                                                                                                                }^{2}(y)}{|x-y|}dxdy\\
\geq & (1-\frac{1}{\eta(a_1,a_2,\beta)} )  [\|(-\Delta)^{\frac{1}{4}}u_{1n}\|_{2}^{2}+\|(-\Delta)^{\frac{1}{4}}u_{2n}\|_{2}^{2} ]+\sum_{i=1}^2\int_{\R ^3}V_{i}(x)|u_{in}|^2 dx,
\end{split}
\end{eqnarray}
which implies that $\{(u_{1n},u_{2n})\}$ is bounded in $\mathcal{X}$ if $\eta(a_1,a_2,\beta) > 1.$ Hence, by Lemma \ref{2:lem1} we may assume that there exists $(u_1,u_2)\in
\mathcal{X}$ such that
\begin{equation*}
\begin{split}
&(u_{1n},u_{2n})\overset{n}{\rightharpoonup}(u_1,u_2) \text{ weakly
in }  \mathcal{X}\,,\\
&(u_{1n},u_{2n}) \overset{n}{\to}(u_1,u_2) \text{ strongly in }
L^q(\R^3)\times L^q(\R^3)\ \text{for } q\in[2,3). \end{split}\end{equation*}
On the other hand, it follows from Lemma 2.2 in \cite{15} that
$$\lim_{n\rightarrow \infty}\int_{\R ^3}\Big(|x|^{-1}\ast |u_{in}|^{2}\Big)u_{in}^{2}dx = \int_{\R ^3}\Big(|x|^{-1}\ast |u_{i}|^{2}\Big)u_{i}^{2}dx, ~~i~ = ~1,2.$$
Note that $\Big(\sqrt{-\Delta + m^{2}}u, u\Big)$ is weakly lower semi-continuous by Lemma A.4 in \cite{5}, then
$$\hat e(a_1,a_2,\beta) \leq E_{a_1,a_2,\beta}(u_1,u_2) \leq \lim_{n\rightarrow \infty}E_{a_1,a_2,\beta}(u_{1n},u_{2n}) = \hat e(a_1,a_2,\beta),$$
since
$\|u_1\|_2^2=\|u_2\|_2^2=1.$ Hence, $E_{a_1,a_2,\beta}(u_1,u_2)=\hat e(a_1,a_2,\beta)$, and $(u_1,u_2)$ is a minimizer of (\ref{eq1.4}).

($\bf ii$) If $\eta(a_1,a_2,\beta) < 1,$ by the definition of (\ref{2:eq1.5}) we may choose $(u_1,u_2)\in\mathcal{M}$ such that $u_1$ and $u_2$ have compact support in $\R^3$ and satisfies
\begin{equation}\label{eq2.2}
\frac{2\|(-\Delta)^{\frac{1}{4}} u_{1}(x) \|_{2}^2+2\|(-\Delta)^{\frac{1}{4}} u_{2}(x) \|_{2}^2}{a_1\int_{\R ^3}\phi_{u_1}u_{1}^{2}dx+a_2\int_{\R ^3}\phi_{u_2}u_{2}^{2}dx+2\beta^{+}\int_{\R ^3}\phi_{u_1}u_{2}^{2}dx}\leq \alpha_0:=\frac{1+\eta(a_1,a_2,\beta)}{2}<1.
\end{equation}
 Let \begin{equation}\label{eq2.3}\bar u_i(x)=\lambda^{\frac{3}{2}} u_i(\lambda x), \text{ for any } \lambda > 0, \ i=1,\,2,\end{equation}
 then $(\bar u_1, \bar u_2)\in\mathcal{M}$. Since $u_i(x)$ has compact support in $\R^3$ and $V_i(x)\in L_{\rm loc}^\infty(\R^3),$ there exists $C > 0$, independent of $\lambda>0$, such that, for $\lambda\to\infty$,
\begin{equation}\label{eq2.4}
\int_{\mathbb{R}^{3}} V_i(x) |\bar u_i|^2dx=\int_{\mathbb{R}^{3}} V_i(\frac{x}{\lambda})|u_i|^2dx\leq C<\infty , \ \, i=1,\,2.
\end{equation}
Moreover, since $\sqrt{-\Delta + m^{2}}  \leq \sqrt{-\Delta} + m$, we infer from (\ref{eq2.2})-(\ref{eq2.4}) that,
\begin{eqnarray}\label{eq2.5}
\begin{split}
E_{a_1,a_2,\beta}(\bar{u}_{1},\bar{u}_{2})&\leq \sum_{i=1}^2\Big(\|(-\Delta)^{\frac{1}{4}}\bar{u}_{i}\|_{2}^{2}+ \int_{\R ^3} V_{i}(x)|\bar{u}_{i}|^2-\frac{a_i}{2}\int_{\R ^3} (|x|^{-1}\ast |\bar{u}_{i}|^{2} )\bar{u}_{i}^{2}dx\Big)\\
& \hspace{2cm} -\beta\iint_{\R ^3\times\R ^3}  \frac{\bar{u}_{1}^{2}(x)\bar{u}_{2}^{2}(y)}{|x-y|}dxdy + 2m\\
&=\lambda\left[\sum_{i=1}^2\Big(\|(-\Delta)^{\frac{1}{4}}u_{i}\|_{2}^{2} -\frac{a_i}{2}\int_{\R ^3} (|x|^{-1}\ast |u_{i}|^{2} )u_{i}^{2}dx\Big) -\beta\int_{\R ^3}\phi_{u_{1}}u_{2}^{2}dx\right]\\
& \hspace{2cm} + \sum_{i=1}^2\int_{\R ^3} V_{i}(x)|\bar{u}_{i}|^2 + 2m\\
&\leq  \lambda(\alpha_0 - 1)\Big(\sum_{i=1}^2\frac{a_i}{2}\int_{\R ^3} (|x|^{-1}\ast |u_{i}|^{2} )u_{i}^{2}dx + \beta\iint_{\R ^3\times\R ^3}  \frac{u_{1}^{2}(x)u_{2}^{2}(y)}{|x-y|}dxdy\Big)\\
& \hspace{2cm} + \sum_{i=1}^2\int_{\R ^3} V_{i}(x)|\bar{u}_{i}|^2 + 2m\\
\rightarrow & -\infty \, \  \text{ as }\,\ \lambda\to +\infty, \ \text{if} \ \beta \geq 0,
\end{split}\end{eqnarray}
 which shows that
$$\hat e(a_1,a_2,\beta)\leq E_{a_1,a_2,\beta}(\bar u_1, \bar u_2)\to-\infty \text{ as } \lambda\to +\infty,$$
and $\hat e(a_1,a_2,\beta)$ can not be attained by a minimizer if $\beta \geq 0$.\\

\indent For $\beta < 0,$ i.e., $\beta^{+} = 0,$ we {\bf claim} that
\[
a_{1} > a^{\ast} \text{ or } a_{2} > a^{\ast} \text{ if } \eta(a_1,a_2,\beta) < 1.
\]
In fact, if $a_{1}, a_{2} \in  (0, a^{\ast} ],$ by Gagliardo-Nirenberg inequality (\ref{GNineq}) we know that, for each $u_i\in H^{\frac{1}{2}}(\R ^3)$ with $\|u_i\|_2^2=1, i = 1,2,$
$$
\aligned J_{a_{1},a_{2},\beta}\big(u_{1},u_{2}\big) & =  \frac{2\|(-\Delta)^{\frac{1}{4}} u_{1}(x) \|_{2}^2+2\|(-\Delta)^{\frac{1}{4}} u_{2}(x) \|_{2}^2}{a_1\int_{\R ^3}\phi_{u_1}u_{1}^{2}dx+a_2\int_{\R ^3}\phi_{u_2}u_{2}^{2}dx}
\geq \frac{a^{\ast}\int_{\R ^3}\phi_{u_1}u_{1}^{2}dx+a^{\ast}\int_{\R ^3}\phi_{u_2}u_{2}^{2}dx}{a_1\int_{\R ^3}\phi_{u_1}u_{1}^{2}dx+a_2\int_{\R ^3}\phi_{u_2}u_{2}^{2}dx}\\
& \geq 1,\endaligned$$
which means that $\eta(a_1,a_2,\beta) \geq 1$, this contradicts our assumption. So,
without loss of generality, we assume that $a_{1} > a^{\ast}$ (the case of $a_{2} > a^{\ast}$ can be done by almost the same way). Take $\varphi \in C_{0}^{\infty}(\mathbb{R}^{3})$ such that
$$\varphi(x) = 1~~ \text{for} ~~|x| \leq \frac{1}{2};~~ \varphi(x) = 0~~ \text{for} ~~|x| > 1; ~~0 \leq \varphi \leq 1.$$
For any $x_{0} \in \mathbb{R}^{3}$ and $R > 1,$ $Q(x)$ is a radially symmetric ground state solution of (\ref{s0000}), let
\begin{equation}\label{eq2.555}
\psi(x)=A_{R}\frac{R^{\frac{3}{2}}}{\|Q\|_{2}}\varphi(x-x_{0})Q \big(R (x-x_{0})\big),
\end{equation}
where $A_{R} > 0$ is chosen such that $\int_{\mathbb{R}^{3}}\psi^{2}dx = 1.$
Then, by (\ref{1:id}) we have
\begin{equation}
\int_{\mathbb{R}^{3}} \psi\sqrt{-\Delta}\psi dx \leq \frac{R}{\|Q\|_{2}^{2}}\Big(\|(-\Delta)^{\frac{1}{4}}Q\|_{2}^{2}\Big) + O(R^{-\frac{5}{2}})= R + O(R^{-\frac{5}{2}}).
 \label{fgytr4}
\end{equation}
\begin{equation}\begin{split}
&\frac{2R}{a^{\ast}} - O(R^{-4}) = \frac{R}{\|Q\|_{2}^{4}}\int_{\R ^3}\Big(|x|^{-1}\ast |Q|^{2}\Big)Q^{2}dx-O(R^{-4}) \leq\iint_{\R ^3\times\R ^3}\frac{\psi^{2}(x)\psi^{2}(y)}{|x-y|}dxdy\\
&\leq \frac{R}{\|Q\|_{2}^{4}}\int_{\R ^3}\Big(|x|^{-1}\ast |Q|^{2}\Big)Q^{2}dx+O(R^{-4})= \frac{2R}{a^{\ast}} + O(R^{-4}).
\label{fgytr456}\end{split}
\end{equation}
By \cite[Lemma 2.2]{nn} we know that there exists $C > 0$ such that
\begin{equation}
|Q(x)| \leq C(1 + |x|)^{-4} \ \mbox{in}\ \ \R^3,
 \label{4:exp}
\end{equation}
and
\begin{equation}
\Big(|x|^{-1}\ast Q^{2}\Big)(x) \leq C(1 + |x|)^{-1} \ \mbox{in}\ \ \R^3. \
 \label{44}
\end{equation}
\indent Since the function $x \mapsto V_{i}(x)\varphi^{2}(\frac{x-x_{0}}{R})$ is bounded and has compact support, by the dominated convergence theorem, we have
\begin{equation}\label{eq2.9991}\begin{split}
\lim_{R\rightarrow\infty}\int_{\mathbb{R}^{3}} V_{i}(x)\psi^{2}(x)dx &= \lim_{R\rightarrow\infty}\int_{\mathbb{R}^{3}} V_{i}(x)A^{2}_{R}\frac{R^{3}}{\|Q\|^{2}_{2}}\varphi^{2}(x-x_{0})Q^{2} \big(R (x-x_{0})\big)dx\\
&= \lim_{R\rightarrow\infty}\frac{A^{2}_{R}}{\|Q\|^{2}_{2}}\int_{\mathbb{R}^{3}}V_{i}(R^{-1}x + x_{0})\varphi^{2}(R^{-1}x)Q^{2}(x)dx\\
&=  V_{i}(x_{0})\ \ \ i = 1,2.
\end{split}\end{equation}
\indent In addition, the Lemma \ref{lemA3} and (\ref{fgytr456}) give that
\begin{equation}\label{bysz}
\int_{\mathbb{R}^{3}}\phi_{\psi}\nu^{2}dx \leq \Big(\int_{\mathbb{R}^{3}}\phi_{\psi}\psi^{2}dx\int_{\mathbb{R}^{3}}\phi_{\nu}\nu^{2}dx\Big)^{\frac{1}{2}}
=\sqrt{\frac{2}{a^{\ast}}}R^{\frac{1}{2}}\Big(\int_{\mathbb{R}^{3}}\phi_{\nu}\nu^{2}dx\Big)^{\frac{1}{2}} + O(R^{-2})
.\end{equation}
\noindent Take $0 \le \nu \in C_0^\infty(\R^3)$ with $\int_{\mathbb{R}^{3}} |\nu|^2dx = 1.$
By $\sqrt{-\Delta + m^{2}}  \leq \sqrt{-\Delta} + m$,  it follows from (\ref{fgytr4})-(\ref{bysz}) that,  for $a_{1} > a^{\ast},$
$$\aligned
E_{a_1,a_2,\beta}(\psi(x),\nu)
& \leq \|(-\Delta)^{\frac{1}{4}}\psi\|_{2}^{2} + \|(-\Delta)^{\frac{1}{4}}\nu\|_{2}^{2} -\frac{a_1}{2}\int_{\mathbb{R}^{3}}\phi_{\psi}\psi^{2}dx - \frac{a_2}{2}\int_{\mathbb{R}^{3}}\phi_{\nu}\nu^{2}dx\\
& \hspace{1.5cm} - \beta\int_{\mathbb{R}^{3}}\phi_{\psi}\nu^{2}dx
+\int_{\R ^3} V_{1}(x)\psi(x)^2dx + \int_{\R ^3} V_{2}(x)\nu^2dx + 2m\\
&\leq R\big(1 - \frac{a_{1}}{a^{\ast}}\big) + \|(-\Delta)^{\frac{1}{4}}\nu\|_{2}^{2} - \frac{a_2}{2}\int_{\mathbb{R}^{3}}\phi_{\nu}\nu^{2}dx
+ |\beta|\sqrt{\frac{2}{a^{\ast}}}R^{\frac{1}{2}}\big(\int_{\mathbb{R}^{3}}\phi_{\nu}\nu^{2}dx\big)^{\frac{1}{2}}\\
&\hspace{1.5cm} +\int_{\R ^3} V_{1}(x)\psi(x)^2dx + \int_{\R ^3} V_{2}(x)\nu^2dx + 2m + O(R^{-\frac{5}{2}}) + O(R^{-2})\\
&\rightarrow -\infty \, \  \text{ as }\,\ R \to +\infty,\endaligned$$
hence,  $\hat e(a_1,a_2,\beta)$ has no minimizer for $\beta <0$, either.

\section{Proofs of the main Theorems}\label{Sec3}

In this section, we come to prove our main Theorems 1.1-1.3.\\

\noindent\textbf{Proof of Theorem \ref{thm1}.}

\indent \textbf{(i):} For any $(u_1,u_2) \in \mathcal{M},$ by Gagliardo-Nirenberg inequality (\ref{GNineq}) and Lemma 2.3, we have
\begin{align*}
J_{a_{1},a_{2},\beta}\big(u_{1},u_{2}\big) &\geq
\frac{a^*\Big(\int_{\R ^3}\phi_{u_1}u_{1}^{2}dx+\int_{\R ^3}\phi_{u_2}u_{2}^{2}dx\Big)}{a_1\int_{\R ^3}\phi_{u_1}u_{1}^{2}dx+a_2\int_{\R ^3}\phi_{u_2}u_{2}^{2}dx+2\beta^{+}\Big(\int_{\R ^3}\phi_{u_1}u_{1}^{2}dx\int_{\R ^3}\phi_{u_2}u_{2}^{2}dx\Big)^{\frac{1}{2}}}\\
&=\frac{a^*\Big(1
+\int_{\R ^3}\phi_{u_2}u_{2}^{2}dx\big/\int_{\R ^3}\phi_{u_1}u_{1}^{2}dx\Big)}{a_1+a_2\int_{\R ^3}\phi_{u_2}u_{2}^{2}dx\big/\int_{\R ^3}\phi_{u_1}u_{1}^{2}dx+2\beta^{+}\big(\int_{\R ^3}\phi_{u_2}u_{2}^{2}dx\big/\int_{\R ^3}\phi_{u_1}u_{1}^{2}dx\big)^\frac{1}{2}}.
\end{align*}
Let
\begin{equation}
\gamma_{a_1,a_2,\beta}(t):=\frac{a^*(1+t^2)}{a_1+a_2t^2+2\beta^{+} t}.
\end{equation}
Then, it is easy to see that,  for any $(u_1,u_2) \in \mathcal{M},$
$$J_{a_{1},a_{2},\beta}\big(u_{1},u_{2}\big) \geq \gamma_{a_1,a_2,\beta}(t_0) ~\text{with}~ t_{0} = \int_{\R ^3}\phi_{u_2}u_{2}^{2}dx\big/\int_{\R ^3}\phi_{u_1}u_{1}^{2}dx.$$
Hence,
 \begin{equation}\label{eq2.77}
 \eta(a_1,a_2,\beta)=\inf\big\{J_{a_{1},a_{2},\beta}\big(u_{1},u_{2}\big): (u_1,u_2) \in \mathcal{M}\big\} \geq \inf_{t\in(0,\infty)} \gamma_{a_1,a_2,\beta}(t).
\end{equation}
By $a_i \in (0, a^*)$ ($i=1,\,2$) and $\beta < \beta_{\ast} = \sqrt {(a^*-a_1)(a^*-a_2)}$, we see that
\[
\gamma_{a_1,a_2,\beta}(t) > 1 \text{ for all } t\in(0,\infty), \text{ and }
\underset{t\to0^+}\lim{\gamma_{a_1,a_2,\beta}}(t)=\frac{a^*}{a_1}>1,\  \underset{t\to\infty}\lim \gamma_{a_1,a_2,\beta}(t)=\frac{a^*}{a_2}>1.
\]
Thus, $\underset{t\in(0,\infty)} \inf \gamma_{a_1,a_2,\beta}(t) > 1$ by the continuity  of $\gamma_{a_1,a_2,\beta}(t),$  then (\ref{eq2.77}) gives that $\eta(a_1,a_2,\beta) > 1.$ So, problem (\ref{eq1.4}) has a minimizer by Theorem \ref{prop:2.1}(i).\\

\textbf{(ii):} We prove now that problem (\ref{eq1.4}) has no minimizer by the following three cases.\\
\indent $\bullet ~\bf{If}$ $a_1 > a^*,$  we take $0 \leq \nu(x) \in C_{0}^\infty \big(\R ^3\big)$ with $\int_{\R ^3} |\nu|^2dx = 1,$ then by the definition of (\ref{5}) as well as (\ref{eq2.555})-(\ref{fgytr456}), we see that, for any $\beta \in \R,$
\begin{equation}\label{equation5}\begin{split}
J_{a_{1},a_{2},\beta}\big(\psi(x),\nu (x) \big)&= \frac{2\|(-\Delta)^{\frac{1}{4}} \psi(x)\|_{2}^2+2\|(-\Delta)^{\frac{1}{4}} \nu(x) \|_{2}^2}{a_1\int_{\R ^3}\phi_{\psi}\psi^{2}dx+a_2\int_{\R ^3}\phi_{\nu}\nu^{2}dx+2\beta^{+}\int_{\R ^3}\phi_{\psi}\nu^{2}dx}\\
&\leq\frac{2R+2\|(-\Delta)^{\frac{1}{4}} \nu(x) \|_{2}^2+ O(R^{-\frac{5}{2}})}{2\frac{a_1}{a^*}R+a_2\int_{\R ^3}\phi_{\nu}\nu^{2}dx- O(R^{-4})} \to\frac{a^*}{a_1} < 1\,\
\text{ as } \,\ R \to +\infty,\end{split}\end{equation}
which gives that $\eta(a_1,a_2,\beta) < 1,$ and then (\ref{eq1.4}) has no minimizer by Theorem \ref{prop:2.1}(ii).\\
\indent $\bullet ~\bf{If}$ $a_2 > a^*$, the proof is almost the same as the above case, we need only change the order of $\psi(x)$ and $\nu(x)$ in \eqref{equation5}.\\
\indent $\bullet ~\bf{If}$ $\beta > \beta^{\ast} = \frac{a^*-a_1}{2} + \frac{a^*-a_2}{2},$ for $\psi(x)$ given by (\ref{eq2.555}), we know from (\ref{fgytr4})-(\ref{fgytr456}) that
\begin{equation}\label{eq2.777}
\begin{split}
&\int_{\mathbb{R}^{3}} \psi\sqrt{-\Delta + m^{2}}\psi dx - \frac{a_{i}}{2}\iint\limits_{\R ^3\times\R ^3}  \frac{\psi^{2}(x)\psi^{2}(y)}{|x-y|}dxdy \\
 \leq  & \frac{R}{\|Q\|_{2}^{2}} (\|(-\Delta)^{\frac{1}{4}}Q\|_{2}^{2}-\frac{a_{i}}{2\|Q\|_{2}^{2}}\iint\limits_{\R ^3\times\R ^3}\frac{Q^{2}(x)Q^{2}(y)}{|x-y|}dxdy )
+ O(R^{-\frac{5}{2}}) + m^{2}O(R^{-1})\\
=&  R\Big(1 - \frac{a_{i}}{\|Q\|_{2}^{2}}\Big) + O(R^{-\frac{5}{2}}) + m^{2}O(R^{-1}) \ \ \ i = 1,2,
\end{split}
\end{equation}
this and (\ref{eq2.9991}) show that
\begin{equation}\label{eq2.741}
E_{a_{1},a_{2},\beta}(\psi, \psi ) \leq R\Big(2 - \frac{a_{1} + a_{2} + 2\beta}{\|Q\|^{2}_{2}}\Big) + \sum_{i=1}^2V_{i}(x_{0}) +O(R^{-\frac{5}{2}}) + m^{2}O(R^{-1}) \overset{R}{\to} -\infty
\end{equation}
if $\beta > \beta^{\ast} = \frac{a^*-a_1}{2} + \frac{a^*-a_2}{2},$ hence problem (\ref{eq1.4}) has no minimizer. \hfill$\square$\\

\noindent \noindent\textbf{Proof of Theorem \ref{thm2}.}\\
\indent Since $0 < a_1 \not= a_2 < a^*,$ $\beta_{\ast} < \beta^{\ast}$ by Remark 1.1 (ii), we know from Lemma \ref{33:lem3.1} that
 \[\eta\big(a_1,a_2,\beta_{\ast}\big) > 1 \text{  and }
 \eta\big(a_1,a_2,\beta^{\ast}\big) \leq 1.
 \]
 However,  $\eta(a_1,a_2,\beta)$ is locally Lipschitz continuous by Lemma \ref{lemA4},  then there exists $\delta := \delta(a_1, a_2) \in \big(\beta_{\ast}, \beta^{\ast}\big]$ such that $\eta\big(a_1,a_2,\beta\big) > 1$ for any $\beta \in \big[\beta_{\ast}, \delta\big),$
and then Theorem \ref{prop:2.1} shows that problem (\ref{eq1.4}) has at least a minimizer for $\beta \in \big[\beta_{\ast}, \delta\big).$ \hfill$\square$\\

Finally, we turn to proving Theorem \ref{thm03}. Before giving its proof, we first establish an estimate for  $\hat e(a^*-\beta,a^*-\beta,\beta)$ based on $V_1(x)$ and $V_2(x)$.
 Under the conditions of Theorem \ref{thm03},  it is easy to see that the set  $\{(a_1,a_2,\beta)\in \mathbb{R}^3: (a_1,a_2,\beta)=(a^*-\beta,a^*-\beta,\beta) \text{ and }   \beta  \in (0, a^*)\}$ forms a segment, and the functional (\ref{f1}) can be rewritten as
\begin{equation}\label{ggg}
\begin{split}
 E_{a_1,a_2,\beta}(u_1, & u_2)
= \sum_{i=1}^2 \int_{\R ^3}  u_i\sqrt{-\Delta + m^{2}}\,u_idx + \sum_{i=1}^2\int_{\R ^3} V_i(x)|u_i|^2dx\\
  &-\frac{a^{\ast}}{2}\int_{\R ^3}\Big(\phi_{u_1}u_{1}^{2}+\phi_{u_2}u_{2}^{2}\Big)dx+\frac{\beta}{2}\int_{\R ^3}\Big(\phi_{u_1}u_{1}^{2}+\phi_{u_2}u_{2}^{2}-2\phi_{u_1}u_{2}^{2}\Big)dx.
\end{split}\end{equation}
 It follows from $\sqrt{-\Delta + m^{2}} \geq \sqrt{-\Delta}$ and (\ref{ggg}) that
\begin{equation}\label{f2}
\begin{split}
E_{a_1,a_2,\beta}(u_1,u_2) &\geq \|(-\Delta)^{\frac{1}{4}}u_1\|^{2} + \|(-\Delta)^{\frac{1}{4}}u_2\|^{2} + \sum_{i=1}^2\int_{\R ^3} V_i(x)|u_i|^2dx\\
&-\frac{a^{\ast}}{2}\int_{\R ^3} (\phi_{u_1}u_{1}^{2}+\phi_{u_2}u_{2}^{2} )dx+\frac{\beta}{2}\int_{\R ^3} (\phi_{u_1}u_{1}^{2}+\phi_{u_2}u_{2}^{2}-2\phi_{u_1}u_{2}^{2} )dx\\
&\geq 0,
\end{split}\end{equation}
since the Gagliardo-Nirenberg inequality (\ref{GNineq}) and Lemma 2.3 imply that
\begin{align*}
\|(-\Delta)^{\frac{1}{4}}u_i\|^{2} - \frac{a^{\ast}}{2}\int_{\R ^3}\phi_{u_i}u_{i}^{2}dx \geq 0, ~~i = 1,2, \nonumber\\
\frac{\beta}{2}\int_{\R ^3}\Big(\phi_{u_1}u_{1}^{2}+\phi_{u_2}u_{2}^{2}-2\phi_{u_1}u_{2}^{2}\Big)dx \geq 0, ~\text{for}~ \beta > 0. \nonumber\end{align*}
\noindent Note that $x_{0} \in \R$ is arbitrary and $\big(a_{1} + a_{2} + 2\beta\big)/\|Q\|_{2}^2 = 2$ in (\ref{eq2.741}), then it follows from (\ref{eq2.741}) and (\ref{f2}) that
\begin{equation}\label{smgy}
0\leq \hat e(a^*-\beta,a^*-\beta,\beta) \leq \inf_{x\in\R^3}\big(V_1(x)+V_2(x)\big).
\end{equation}

\noindent{\textbf{Proof of Theorem \ref{thm03}.}}

\indent \textbf{(i):} Since $\inf_{x\in \R^3}\big(V_1(x)+V_2(x)\big) = 0$ and $V_i(x) \geq 0$ $(i = 1,2),$ we know that there exists  $x_0 \in \R^3$ such that  $V_1(x_0) = V_2(x_0) = 0$. It follows from (\ref{smgy}) that
\begin{equation}\label{5:ea*}
0\le \hat e(a^*-\beta,a^*-\beta,\beta)\leq \inf_{x\in\R^3}\big(V_1(x)+V_2(x)\big) = V_1( x_0)+V_2( x_0) = 0,
\end{equation}
i.e.,
$$\hat e(a^*-\beta,a^*-\beta,\beta)=0.$$
\indent If $\hat e(a^*-\beta,a^*-\beta,\beta) = 0$ has a minimizer $(\hat u_1,\hat u_2) \in \mathcal{M},$ then, it follows from (\ref{f2}) that
\begin{equation}\label{5mmnn}
\int_{\R ^3}\phi_{\hat u_1}\hat u_{1}^{2} = \int_{\R ^3}\phi_{\hat u_2}\hat u_{2}^{2} ~~\text{and}~~\int_{\R ^3}\phi_{\hat u_1}\hat u_{2}^{2} = \big(\int_{\R ^3}\phi_{\hat u_1}\hat u_{1}^{2}\int_{\R ^3}\phi_{\hat u_2}\hat u_{2}^{2}\big)^{\frac{1}{2}}
,\end{equation}
\begin{equation}\label{gbv}
\|(-\Delta)^{\frac{1}{4}}\hat u_1\|_{2}^{2}=\frac{a^*}{2}\iint_{\R ^3\times\R ^3}  \frac{\hat u_1^{2}(x)\hat u_1^{2}(y)}{|x-y|}dxdy\,\ {\rm
and}\,\ \int_{\R^3}V_1(x)|\hat u_1|^2dx=0.
\end{equation}
Since $\hat u_i$ can be assumed to be nonnegative,  it follows from Lemma 2.3 and (\ref{5mmnn}) that $\hat u_1(x) \equiv \hat u_2(x)$ in $\R^3,$ which however contradicts (\ref{gbv}) since the first equality of (\ref{gbv}) shows that $\hat u_1(x)=Q(x)>0$ (up to a translation) , but we know from the second equality of (\ref{gbv})
that $\hat u_1(x)$ is compactly supported in $\R^3$.\\

\textbf{(ii):} Let $\mathcal{M}$ and $\mathcal{X}$ be given by (\ref{eq1.4}), it is clear that $(\mathcal{M},d)$ is a complete metric space with
$$d(\vec u,\vec v):=\|\vec u-\vec v\|_\mathcal{X}, \quad   \vec u,\, \vec v\in \mathcal{M}\,$$
where $\vec u = (u_1, u_2) \in \mathcal{X} , \vec v = (v_1, v_2) \in \mathcal{X},$ and $\|\vec u\|_\mathcal{X}=\big(\|u_1\|_{\mathcal{H}_1}^2+\|u_2\|_{\mathcal{H}_2}^2\big)^\frac{1}{2}.$
By the Ekeland's variational principle \cite[Theorem  5.1]{13}, there is a minimizing sequence
 $\{\vec u_n=(u_{1n},u_{2n})\}\subset \mathcal{M}$ for $\hat e(a^*-\beta,a^*-\beta,\beta)$ such that
\begin{eqnarray}\label{eqE1}
&&\hat e(a^*-\beta,a^*-\beta,\beta)\leq E_{a^*-\beta,a^*-\beta,\beta}(\vec u_n)\leq \hat e(a^*-\beta,a^*-\beta,\beta)+\frac{1}{n}, \label{mn}\\
&&E_{a^*-\beta,a^*-\beta,\beta}(\vec v)\geq E_{a^*-\beta,a^*-\beta,\beta}(\vec u_n)-\frac{1}{n}\|\vec u_n-\vec v\|_\mathcal{X} \quad \text{for}\quad \vec v \in \mathcal{M}.\label{eqE2}
\end{eqnarray}
Since the compact embedding of Lemma \ref{2:lem1}, to prove that $\hat e(a^*-\beta,a^*-\beta,\beta)$ is attained, we need only to show that $\{\vec u_n=(u_{1n},u_{2n})\}$ is bounded in $\mathcal{X}.$\\
\indent In fact, if $\{\vec u_n=(u_{1n},u_{2n})\}$ is unbounded in $\mathcal{X}$, then, passing to a subsequence if necessary, we may assume that $\|\vec u_n\|_\mathcal{X} \xrightarrow{n} +\infty$.
 Applying Gagliardo-Nirenberg inequality (\ref{GNineq}) and Lemma 2.3, it follows from (\ref{f2}) and (\ref{mn}) that
\begin{equation}\label{5:vi}
\sum_{i=1}^2\int_{\mathbb{R}^{3}} V_i(x)|u_{in}|^2dx\leq E_{a^*-\beta,a^*-\beta,\beta}(\vec u_n)\leq \hat e(a^*-\beta,a^*-\beta,\beta)+\frac{1}{n}.
\end{equation}
Then, by the definition (\ref{eq1.789}) and $\|\vec u_n\|_\mathcal{X} \xrightarrow{n} +\infty$  we know that
\begin{equation}\label{eqE4}\|(-\Delta)^{\frac{1}{4}}u_{1n}\|_{2}^{2}+\|(-\Delta)^{\frac{1}{4}}u_{2n}\|_{2}^{2}\xrightarrow{n}\infty ,\end{equation}
and then it follows from $\sqrt{-\Delta + m^{2}} \geq \sqrt{-\Delta}$ that
\begin{equation}\label{xblp}
\int_{\R ^3} u_{1n}\sqrt{-\Delta + m^{2}}u_{1n}dx + \int_{\R ^3} u_{2n}\sqrt{-\Delta + m^{2}}u_{2n}dx \xrightarrow{n}\infty.
\end{equation}

Now, we claim that, for $i = 1, 2,$
\begin{eqnarray}\label{5:com1}
&&\int_{\R ^3} u_{in}\sqrt{-\Delta + m^{2}}u_{in}dx \approx \frac{a^*}{2}\iint_{\R ^3\times\R ^3}  \frac{ u_{in}^{2}(x) u_{in}^{2}(y)}{|x-y|}dxdy \xrightarrow{n} +\infty,\\
&&\iint_{\R ^3\times\R ^3}  \frac{ u_{1n}^{2}(x) u_{1n}^{2}(y)}{|x-y|}dxdy \approx\iint_{\R ^3\times\R ^3} \frac{ u_{2n}^{2}(x) u_{2n}^{2}(y)}{|x-y|}dxdy,\label{5:com2}\\
&&\iint_{\R ^3\times\R ^3}  \frac{ u_{1n}^{2}(x) u_{2n}^{2}(y)}{|x-y|}dxdy \approx\iint_{\R ^3\times\R ^3} \frac{ u_{2n}^{2}(x) u_{2n}^{2}(y)}{|x-y|}dxdy,\label{5:com3}
\end{eqnarray}
here we denote by $f_n\approx g_n$ two function sequences satisfying $f_n/g_n\to1$ as $n\to\infty$.\\
\indent In fact, since (\ref{xblp}), without loss of generality, we suppose that
$$\int_{\R ^3} u_{1n}\sqrt{-\Delta + m^{2}}u_{1n}dx \xrightarrow{n} \infty.$$
By (\ref{eqE1}) and Lemma 2.3, we know that
\begin{equation}\label{5:ui}
0\leq \int_{\R ^3} u_{1n}\sqrt{-\Delta + m^{2}}\,u_{1n}dx-\frac{a^*}{2}\iint\limits_{\R ^3\times\R ^3}  \frac{ u_{1n}^{2}(x) u_{1n}^{2}(y)}{|x-y|}dxdy\leq \hat e(a^*-\beta,a^*-\beta,\beta)+\frac{1}{n},
\end{equation}
with $e(a^*-\beta,a^*-\beta,\beta)+\frac{1}{n} < C,$ for some constant $C > 0$ and $n$ large enough.\\
\noindent Multiplying (\ref{5:ui}) by
$1 \Big{/} \int_{\R ^3} u_{1n}\sqrt{-\Delta + m^{2}}\,u_{1n}dx,$ we have, for $n$ large enough,
\begin{equation}\label{5:u1cxz}\begin{split}
0 &\leq 1 - \frac{a^*}{2}\iint_{\R ^3\times\R ^3}  \frac{ u_{1n}^{2}(x) u_{1n}^{2}(y)}{|x-y|}dxdy\bigg/\int_{\R ^3} u_{1n}\sqrt{-\Delta + m^{2}}\,u_{1n}dx\\
&\leq C \Big{/} \int_{\R ^3} u_{1n}\sqrt{-\Delta + m^{2}}\,u_{1n}dx \rightarrow 0, ~\text{as}~ n \rightarrow +\infty.\end{split}\end{equation}
This means that
\begin{equation}\label{5:u1}\begin{split}
&\frac{a^*}{2}\iint_{\R ^3\times\R ^3}  \frac{ u_{1n}^{2}(x) u_{1n}^{2}(y)}{|x-y|}dxdy = \frac{a^*}{2} \int_{\R^3}\phi_{u_{1n}}u_{1n}^{2}dx  \xrightarrow{n}  \infty, \text{ and }\\
&\int_{\R ^3} u_{1n}\sqrt{-\Delta + m^{2}}u_{1n}dx\bigg/\{\frac{a^*}{2}\iint_{\R ^3\times\R ^3}  \frac{ u_{1n}^{2}(x) u_{1n}^{2}(y)}{|x-y|}dxdy\} \xrightarrow{n} 1.\end{split}\end{equation}
Applying (\ref{eqE1}), Gagliardo-Nirenberg inequality (\ref{GNineq}) and Lemma \ref{lemA3}, we see that
\begin{equation}\label{5:eq}
0\leq\frac{\beta}{2}\int_{\R ^3} (\phi_{u_{1n}}u_{1n}^{2}+\phi_{u_{2n}}u_{2n}^{2}-2\phi_{u_{1n}}u_{2n}^{2} )dx \leq \hat e(a^*-\beta,a^*-\beta,\beta)+\frac{1}{n}, ~\text{for}~\beta > 0.
\end{equation}
Then, for $n$ large,
\begin{equation}\label{5:u15}\begin{split}
\frac{\beta}{2}\bigg| (\int_{\R^3}\phi_{u_{1n}}u_{1n}^{2}  dx)^{\frac{1}{2}}- (\int_{\R^3}\phi_{u_{2n}}u_{2n}^{2} dx)^{\frac{1}{2}}\bigg|^{2}
&\leq \frac{\beta}{2}\int_{\R ^3} (\phi_{u_{1n}}u_{1n}^{2}+\phi_{u_{2n}}u_{2n}^{2}-2\phi_{u_{1n}}u_{2n}^{2} )dx\\
& \leq \hat e(a^*-\beta,a^*-\beta,\beta)+\frac{1}{n} \leq C,
\end{split}\end{equation}
and by $\int_{\R^3}\phi_{u_{1n}}u_{1n}^{2} dx \xrightarrow{n} \infty,$ we know that
\begin{equation}\label{5:u12}
\int_{\R^3}\phi_{u_{2n}}u_{2n}^{2}\xrightarrow{n}\infty \quad \text{and}\,\ \int_{\R^3}\phi_{u_{1n}}u_{1n}^{2}\bigg/\int_{\R^3}\phi_{u_{2n}}u_{2n}^{2}\xrightarrow{n}1.
\end{equation}
Similar to (\ref{5:ui}), we have also
$$0\leq \int_{\R ^3} u_{2n}\sqrt{-\Delta + m^{2}}u_{2n}dx-\frac{a^*}{2}\int_{\R^3}\phi_{u_{2n}}u_{2n}^{2} \leq \hat e(a^*-\beta,a^*-\beta,\beta)+\frac{1}{n},$$
then, by (\ref{5:u12}) we get that
\begin{equation}\label{5:u2} \int_{\R ^3} u_{2n}\sqrt{-\Delta + m^{2}}\, u_{2n}dx \xrightarrow{n} \infty  \text{ and }  \int_{\R ^3} u_{2n}\sqrt{-\Delta + m^{2}}\, u_{2n}dx\bigg/\{\frac{a^{\ast}}{2}\int_{\R^3}\phi_{u_{2n}}u_{2n}^{2} dx \} \xrightarrow{n}1. \end{equation}
Hence,  (\ref{5:com1}) and (\ref{5:com2}) follows from (\ref{5:u1}), (\ref{5:u2}) and (\ref{5:u12}).
So, to complete the proof of our claim (\ref{5:com1})-(\ref{5:com3}), we still need to prove (\ref{5:com3}).\\
\indent Since (\ref{5:eq}) and $\beta > 0$, we have, for some constant $C>0$,
$$0 \leq \int_{\R ^3}\Big(\phi_{u_{1n}}u_{1n}^{2}+\phi_{u_{2n}}u_{2n}^{2}-2\phi_{u_{1n}}u_{2n}^{2}\Big) dx \leq C,$$
that is,
\begin{equation}\label{sbbs}0 \leq 1 + \frac{\int_{\R ^3}\phi_{u_{1n}}u_{1n}^{2}dx}{\int_{\R ^3}\phi_{u_{2n}}u_{2n}^{2}dx} - 2\frac{\int_{\R ^3}\phi_{u_{1n}}u_{2n}^{2}dx}{\int_{\R ^3}\phi_{u_{2n}}u_{2n}^{2}dx} \leq \frac{C}{\int_{\R ^3}\phi_{u_{2n}}u_{2n}^{2}dx} \xrightarrow{n} 0.\end{equation}
This and (\ref{5:u12}) imply that
$$\frac{\int_{\R ^3}\phi_{u_{1n}}u_{2n}^{2}}{\int_{\R ^3}\phi_{u_{2n}}u_{2n}^{2}} \xrightarrow{n} 1,$$ which gives (\ref{5:com3}). Hence, our claim (\ref{5:com1})-(\ref{5:com3}) are proved.

As we mentioned in the beginning of the proof, in order to finish the proof of part \textbf{(ii)} we need only to get a contradiction under the condition (\ref{xblp}). For this purpose,  let
$$\epsilon_n^{-1}:=\int_{\R ^3} u_{1n}\sqrt{-\Delta + m^{2}} \, u_{1n}dx \text{ and } \tilde{w}_{in}(x) = \epsilon_n^{\frac{3}{2}}u_{in}(\epsilon_n x),~i = 1,2.$$
Then,
 \[
 \epsilon_n \xrightarrow{n} 0 \text{  by (\ref{5:com1}),  and }  \|\tilde{w}_{in}\|_{2} = \|u_{in}\|_{2} = 1~(i = 1,2).
  \]
  Moreover, we claim that there exist $R_0 > 0,$ $\eta > 0$ and a sequence $\{y_{\epsilon_n}\}\subset\R^3$ such that
\begin{equation}\label{5:scal}
\liminf_{n\rightarrow \infty}\int_{B_{R_0}(y_{\epsilon_n})}|\tilde{w}_{1n}(x)|^2dx \geq \eta > 0.
\end{equation}
In fact, if (\ref{5:scal}) fails, then, passing to a subsequence if necessary, we have
\begin{equation}\label{5:s11}
\lim_{n\rightarrow \infty}\sup_{y \in \mathbb{R}^3}\int_{B_{R}(y)} |\tilde{w}_{1n}(x)|^2dx = 0, \text{ for any } R > 0.
\end{equation}
Then, it follows from the vanishing lemma due to P.L.Lions \cite[Lemma I.1]{12} that $\tilde{w}_{1n} \xrightarrow{n} 0$ strongly in $L^{r}(\mathbb{R}^{3})$ with $r \in (2, 3),$ which then gives that
\begin{equation}\label{5555s}
\int_{\R^3}\phi_{\tilde{w}_{1n}}\tilde{w}_{1n}^{2} \rightarrow 0 ~\text{as}~ n \rightarrow \infty
\end{equation}
by using Lemma \ref{lemA2} with $p = q = \frac{6}{5}$ and $t = 1.$ However, by (\ref{5:com1}) and the definition of $\epsilon_{n},$ we see that
\begin{equation}\label{5:s12}
\int_{\R^3}\phi_{\tilde{w}_{1n}}\tilde{w}_{1n}^{2} = \epsilon_{n}\int_{\R^3}\phi_{u_{1n}}u_{1n}^{2} \xrightarrow{n} \frac{2}{a^{\ast}} \neq 0,
\end{equation}
which contradicts (\ref{5555s}). So, (\ref{5:scal}) is proved.

{\it We should mention that,  it is easy to show that  (\ref{5:scal}) also holds for  $\{ \tilde{w}_{2n} \}$ by almost the same procedures as above, but it is difficult to know whether we can get the same sequence $\{y_{\epsilon_n}\}$ in both cases, which is crucial for our problem.  Unfortunately,  the approaches used in papers \cite{8,9}  do not work in our case}. To overcome this difficulty, some new strategies are used in this paper.

For the same $\{y_{\epsilon_n}\}$ as in (\ref{5:scal}), we define now
$$w_{in}(x) := \tilde{w}_{in}\big(x + y_{\epsilon_n}\big) = \epsilon_n^{\frac{3}{2}}u_{in}\big(\epsilon_n x + \epsilon_n y_{\epsilon_n}\big),~ i = 1, 2.$$
Then, we have
\begin{equation}\label{5:ssss}
\liminf_{n\rightarrow \infty}\int_{B_{R_0}(0)}|w_{1n}(x)|^2dx \geq \eta > 0.
\end{equation}
Furthermore, $\{\epsilon_ny_{\epsilon_n}\}$ is a bounded sequence in $\R^3.$ Otherwise, by (\ref{5:vi}) we know that
\begin{equation}\label{5:bound}
\begin{split}
\int_{\mathbb{R}^{3}} V_1(x)|u_{1n}|^2dx & = \int_{\mathbb{R}^{3}} V_1(\epsilon_nx+\epsilon_ny_{\epsilon_n})|w_{1n}|^2dx
\leq \sum_{i=1}^2\int_{\R ^3} V_i(x)|u_{in}|^2dx \\
&\leq \hat e(a^*-\beta,a^*-\beta,\beta)+\frac{1}{n} \leq C, \text{ for some } C>0,
\end{split}
\end{equation}
which leads to a contradiction by Fatou's lemma and (\ref{5:ssss}) as well as $V_1(x) \xrightarrow{|x|\to\infty} +\infty.$\\

For any $\varphi(x)\in C_0^\infty(\R^3)$, define
\begin{equation*}\tilde \varphi(x)=\varphi\big(\frac{x-\epsilon_ny_{\epsilon_n}}{\epsilon_n}\big),\,\
j(\tau,\sigma)=\frac{1}{2}\int_{\R^3}\big|u_{1n}+\tau u_{1n}+\sigma\tilde \varphi\big|^2dx.
\end{equation*}
Then, $j(\tau,\sigma)$ satisfies
$$j(0,0)=\frac{1}{2},\,\ \frac{\partial j(0,0)}{\partial \tau}= \int_{\R^3}|u_{1n}|^2dx=1\,\ \text{and}\,\ \frac{\partial j(0,0)}{\partial \sigma}=\int_{\R^3} u_{1n}\tilde \varphi dx.$$
By the implicit function theorem, there exist a constant $\delta_n > 0$ and a function $\tau(\sigma)\in C^1\big((-\delta_n,\delta_n), \R\big)$ such that
$$\tau(0)=0,\,\ \tau'(0)=-\int_{\R^3} u_{1n}\tilde \varphi dx, \,\ \text{and}\ \ j(\tau(\sigma),\sigma)=j(0,0)=\frac{1}{2},$$
then
$$\big(u_{1n}+\tau(\sigma) u_{1n}+\sigma\tilde \varphi, u_{2n}\big)\in \mathcal{M},  \text{ for }  \sigma\in (-\delta_n,\delta_n).$$
Applying (\ref{eqE2}), it is clear that
$$E_{a^*-\beta,a^*-\beta,\beta}(u_{1n}+\tau(\sigma) u_{1n}+\sigma\tilde \varphi, u_{2n})-E_{a^*-\beta,a^*-\beta,\beta}(u_{1n}, u_{2n})\geq-\frac{1}{n}\|(\tau(\sigma)u_{1n}+\sigma\tilde\varphi,0)\|_\mathcal{X}.$$
Using this fact and letting $\sigma\to 0^+$ and $\sigma\to 0^-,$ respectively, we then have
\begin{equation}\label{5:de}\Big|\big\langle E_{a^*-\beta,a^*-\beta,\beta}'(u_{1n},u_{2n}), (\tau'(0)u_{1n}+\tilde \varphi, 0)\big\rangle\Big|\leq\frac{1}{n}\|\tau'(0)u_{1n}+\tilde \varphi\|_{\mathcal{H}_1}.
\end{equation}
By the definitions of $\epsilon_n$ and $w_{1n},$ it is not difficult to see that
\begin{equation}\label{5:de2}
\tau'(0)=-\int_{\R^3} u_{1n}\tilde\varphi dx=-\epsilon_n^\frac{3}{2}\int_{\R^3} w_{1n}\varphi dx\,,\quad \|\tau'(0)u_{1n}+\tilde \varphi\|_{\mathcal{H}_1} \leq C\epsilon_n.
\end{equation}
Denote
\begin{align}\label{5:de02}
\mu_{1n}:=   \langle E_{a^*-\beta,a^*-\beta,\beta}'(u_{1n},u_{2n}),  & (u_{1n}, 0)\rangle
= 2\int_{\R ^3} u_{1n}\sqrt{-\Delta + m^{2}} \, u_{1n}dx + 2\int_{\R^3} V_1(x)u_{1n}^2 dx \nonumber \\
&     - 2a_{1}\int_{\R^3}\phi_{u_{1n}}u_{1n}^2 dx - 2\beta\int_{\R^3}\phi_{u_{2n}}u_{1n}^2 dx,
\end{align}
Then, for $a_{1} = a^*-\beta$, it follows from (\ref{5:com1})-(\ref{5:com3}) and (\ref{5:bound}) that
\begin{equation}\label{5:de002}
\mu_{1n} \approx -a^{\ast}\int_{\R^3}\phi_{u_{1n}}u_{1n}^2 + 2\int_{\R^3} V_1(x)u_{1n}^2 , \text{ and }  \mu_{1n}\epsilon_n \xrightarrow{n} -2.
\end{equation}
Therefore, using (\ref{5:de}), (\ref{5:de2})-(\ref{5:de002}), we see that
\begin{align}\label{5:de3}
&\Big|\big\langle E_{a^*-\beta,a^*-\beta,\beta}'(u_{1n},u_{2n}), (\tau'(0)u_{1n}+\tilde \varphi, 0)\big\rangle\Big|
=2\sqrt{\epsilon_n}\bigg|\int_{\R^3}\varphi\sqrt{-\Delta + \epsilon_n^{2}m^{2}}w_{1n}dx \nonumber \\
& +\epsilon_n \int_{\R^3} V_1(\epsilon_n x+\epsilon_ny_{\epsilon_n})w_{1n}\varphi dx-\frac{\mu_{1n}\epsilon_n}{2}\int_{\R^3}w_{1n}\varphi dx
-\int_{\R^3}(a_1 \phi_{w_{1n}}w_{1n}+ \beta \phi_{w_{2n}}w_{1n}) \varphi dx\bigg| \nonumber\\
&\hspace{5.5cm}  \leq  \frac{1}{n}\|\tau'(0)u_{1n}+\tilde \varphi\|_{\mathcal{H}_1} \leq \frac{C\epsilon_n}{n} \rightarrow 0 \ \ \text{as} \ n \rightarrow \infty,
\end{align}
which then implies that, for any $\varphi(x)\in C_0^\infty(\R^3),$
\begin{align}\label{6.dede}
&\int_{\R^3}\varphi\sqrt{-\Delta + \epsilon_n^{2}m^{2}} \, w_{1n}dx +\epsilon_n\int_{\R^3} V_1(\epsilon_n x+\epsilon_ny_{\epsilon_n})w_{1n}\varphi dx-\frac{\mu_{1n}\epsilon_n}{2}\int_{\R^3}w_{1n}\varphi dx \nonumber \\
& \hspace{5cm}  -a_1\int_{\R^3}\phi_{w_{1n}}w_{1n}\varphi - \beta\int_{\R^3}\phi_{w_{2n}}w_{1n}\varphi dx \xrightarrow{n} 0.
\end{align}
Since $\{w_{in}\}$ $(i = 1,2)$ are bounded in $H^\frac{1}{2}(\R^3),$ we may assume that
$$w_{in}\overset{n}{\rightharpoonup} w_i~ \text{weakly~ in} ~H^\frac{1}{2}(\R^3),~\text{for ~some}~ w_i \in H^\frac{1}{2}(\R^3),~i = 1,2$$
and $w_1 \not\equiv 0$ by (\ref{5:ssss}).  In what follows, we come to prove that $w_2 \not\equiv 0.$

In fact, by (\ref{5:de002}) and (\ref{6.dede}), we know that $w_1$ and $w_{2}$ satisfy, in the weak sense, the following equation
$$\sqrt{-\Delta}w_{1} + w_{1} = a_1\phi_{w_{1}}w_{1} + \beta\phi_{w_{2}}w_{1}, ~a_1 = a^{\ast} - \beta.$$
Similarly, by noting $a_2 = a_1,$ $w_1$ and $w_{2}$ also satisfy
$$\sqrt{-\Delta}w_{2} + w_{2} = a_1\phi_{w_{2}}w_{2} + \beta\phi_{w_{1}}w_{2}.$$
Hence, $\big(w_1, w_2\big)$ is a weak solution for the system
\begin{equation}\label{5:lim}
\begin{cases}
 \sqrt{-\Delta}u + u = a_1\phi_{u}u + \beta\phi_{v}u
, \\
 \sqrt{-\Delta}v + v = a_1\phi_{v}v + \beta\phi_{u}v.
\end{cases}
\end{equation}
By contradiction, if $w_2 \equiv 0,$ then $w_{1}$ satisfies the following equation
\begin{equation}\label{drthjui}
\sqrt{-\Delta}w_{1} + w_{1} = a_1\phi_{w_{1}}w_{1}.\end{equation}
Applying Theorem 6.1 in \cite{bc},  $w_{1}$ satisfies the following type Pohozaev  identity
$$\int_{\R^3} |(-\Delta)^{\frac{1}{4}}w_1|^{2}dx + \frac{3}{2}\int_{\R^3} w_1^2 dx =  \frac{5}{4}a_1\iint_{\R ^3\times\R ^3} \frac{w_1^2(x)w_1^2(y)}{|x - y|}dxdy.$$
And,  by  (\ref{drthjui})  we have
$$\int_{\R^3} |(-\Delta)^{\frac{1}{4}}w_1|^{2}dx + \int_{\R^3} w_1^2 dx = a_1\iint_{\R ^3\times\R ^3} \frac{w_1^2(x)w_1^2(y)}{|x - y|}dxdy\big),$$
Then,
\begin{equation}\label{dyrvl}
\int_{\R^3} |(-\Delta)^{\frac{1}{4}}w_1|^{2}dx = \int_{\R^3} w_1^2 dx = \frac{a_1}{2}\iint_{\R ^3\times\R ^3} \frac{w_1^2(x)w_1^2(y)}{|x - y|}dxdy.\end{equation}
This and  Gagliardo-Nirenberg inequality (\ref{GNineq}) show that
$$\frac{2}{a_1}\|w_1\|^{2}_{2} =\iint_{\R ^3\times\R ^3} \frac{w_1^2(x)w_1^2(y)}{|x - y|}dxdy \leq \frac{2}{a^\ast} \|(-\Delta)^{\frac{1}{4}} w_1\|_{2}^2\|w_1\|^{2}_{2} = \frac{2}{a^\ast}\|w_1\|^{4}_{2},$$
hence, we get from $a_1 = a^\ast - \beta \in \big(0, a^\ast\big)$  that
$$\|w_1\|^{2}_{2} \geq \frac{a^\ast}{a_1} > 1,$$
which is impossible since $\|w_1\|^{2}_{2} \leq \lim\limits_{n\rightarrow \infty}\|w_{1n}\|^{2}_{2} = 1.$ So, $w_2 \not\equiv 0$,  that is,  $\big(w_1, w_2\big)$ is a nontrivial solution of the system (\ref{5:lim}).

 Let $\mathcal{X}_0=H^{\frac{1}{2}}(\R ^3) \times H^{\frac{1}{2}}(\R ^3)$,  $\|u\|_{H^\frac{1}{2}}^2=\int_{\R^3} (|(-\Delta)^{\frac{1}{4}}u|^{2} + u^{2} )dx $ and
\begin{align*}
\mathcal{N} = \{(u, v) \in \mathcal{X}_0:  u,v \neq 0, &\|u\|_{H^\frac{1}{2}}^2 =  \int_{\R^3} (a_1\phi_{u}u^{2} + \beta\phi_{v}u^{2} )dx,
\|v\|_{H^\frac{1}{2}}^2 =  \int_{\R^3} (a_1\phi_{v}v^{2} + \beta\phi_{u}v^{2} )dx \}
\end{align*}
and
\begin{equation}\label{5:limbvc}
J(u, v):= \frac{1}{2}(\|u\|_{H^\frac{1}{2}}^2 + \|v\|_{H^\frac{1}{2}}^2)
-\frac{1}{4}\int_{\R^3}\Big(a_1\phi_{u}u^{2}+ a_1\phi_{v}v^{2}+2\beta\phi_{u}v^{2}\Big)\end{equation}
Then, for $Q > 0$ being a radially symmetric ground state solution of (\ref{s0000}), it follows from Lemmas 4.1 and 4.2 in the Appendix that
\begin{equation}\label{ghbc}
 (u_{0}, v_{0} ) = \begin{cases}\begin{split}
 (\frac{1}{\sqrt{a^{\ast}}}Q(x), \frac{1}{\sqrt{a^{\ast}}}Q(x) ), ~~~&a_1 \neq \beta \text{ and }  ~a_1 + \beta = a^\ast,\\
 (\frac{1}{\sqrt{a_1}}Q\sin \theta, \frac{1}{\sqrt{a_1}}Q\cos \theta ), ~~~&\theta \in (0, 2\pi) \text{ and } a_1 = \beta = \frac{a^{\ast}}{2},
\end{split}\end{cases}\end{equation}
is a minimizer for the following minimization problem
$$\inf \{J(u, v): (u,v)\in \mathcal{N}\}.$$
Using Theorem 6.1 of \cite{bc}, we know that a solution $(u,v)$ of system (\ref{5:lim}) satisfies the following type Pohozaev identity
\begin{equation}\label{pohozaev}
\int_{\R^3}  |(-\Delta)^{\frac{1}{4}}u|^{2}+|(-\Delta)^{\frac{1}{4}}v|^{2}   dx + \frac{3}{2}\int_{\R^3}  u^{2} + v^{2}  dx =
\frac{5}{4}\int_{\R^3} (a_1\phi_{u}u^{2}+ a_1\phi_{v}v^{2}+2\beta\phi_{u}v^{2} )dx.
\end{equation}
This together with  $(u,v) \in  \mathcal{N}$, we have
\begin{equation}\label{equation1}
\int_{\R^3} |(-\Delta)^{\frac{1}{4}}u|^{2}+|(-\Delta)^{\frac{1}{4}}v|^{2} dx  = \int_{\R^3} u^{2} + v^{2} dx
=\frac{1}{2}\int_{\R^3} a_1\phi_{u}u^{2}+ a_1\phi_{v}v^{2}+2\beta\phi_{u}v^{2} dx.
\end{equation}
Therefore, it follows from (\ref{5:limbvc})-(\ref{ghbc}) and (\ref{equation1}) that
\begin{equation}\label{hsw}
J(w_{1}, w_{2}) = \frac{1}{2}\int_{\R^3}\Big(w_{1}^{2} + w_{2}^{2}\Big) \geq J(u_{0},v_{0}) = \frac{1}{2}\int_{\R^3}\Big(u_{0}^{2} + v_{0}^{2}\Big) = 1.\end{equation}
Note that $\|w_{in}\|_2^2 \equiv 1$ and $w_i \not \equiv 0$, we know that $0< \|w_{i}\|_2^2 \leq 1,$ and it then follows from (\ref{hsw}) that
\begin{equation}\label{5:cons}
\|w_{i}\|_2^2 = \int_{\R^3}w_{i}^2 = 1,\,\ \mbox{where}\,\ i=1,\, 2.
\end{equation}
\noindent Thus,
\begin{equation}\label{5:strong}w_{in}\xrightarrow{n} w_i\,\ \text{strongly in}\,\ L^2(\R^3), \,\ i=1,\, 2.\end{equation}
Since $\{\epsilon_n y_{\epsilon_n}\}$ is bounded in $\R^3$, 
we may assume that $\epsilon_n y_{\epsilon_n}\xrightarrow{n}z_0\in \R^3$.
Applying Fatou's lemma, it follows from (\ref{5:bound}), (\ref{5:cons}) and (\ref{5:strong})  that
\begin{equation*}
\begin{split}
\hat{e}(a^*-\beta,a^*-\beta,\beta)&\geq \sum_{i=1}^2\int_{\R^3}\lim_{n\to\infty} V_i(\epsilon_nx+\epsilon_ny_{\epsilon_n})|w_{in}|^2dx\\
&= \sum_{i=1}^2\int_{\R^3} V_i(z_0)|w_{i}|^2dx=V_1(z_0)+V_2(z_0),
\end{split}
\end{equation*}
which contradicts (\ref{eq1.1400003}). Hence, $\{\big(u_{1n},u_{2n}\big)\}$ is bounded in $\mathcal{X}$ and the proof of Theorem \ref{thm03} (ii) is completed. ~~\hfill$\square$

\section{Appendix}\label{Sec4}

\indent For the sake of completeness, in this appendix we come to prove two lemmas which have been used in the proof of Theorem \ref{thm03} (ii). As in the Introduction, let $Q(x) > 0$ be a ground state solution of equation (\ref{s0000}) and $a^\ast = \|Q\|_{2}^2.$ We consider now the following system
\begin{equation}\label{nb666}
\begin{cases}
 \sqrt{-\Delta}u + u = a\phi_{u}u + \beta\phi_{v}u ~~~~(u, v) \in \mathcal{X}_0 = H^{\frac{1}{2}}(\R ^3) \times H^{\frac{1}{2}}(\R ^3),\\
 \sqrt{-\Delta}v + v = a\phi_{v}v + \beta\phi_{u}v ~~~~a,~\beta \in \mathbb{R}^+,
\end{cases}
\end{equation}
where $\phi_{f}(x)$ is defined as (\ref{5}). Our following two lemmas show that certain type ground state solutions for the system (\ref{nb666}) can be constructed by virtue of $Q(x)$. If the fractional operator $\sqrt{-\Delta}$ is replaced by the Laplace operator $-\Delta$ in (\ref{s0000}) and (\ref{nb666}), it was proved in paper \cite[Theorem 1.3]{1} that $ (\sqrt{k}Q, \sqrt{l}Q )$ is a ground state solution to (\ref{nb666}) (with $-\Delta$ operator) for some $k > 0, l > 0$ and $a \neq \beta.$ Motivated by the proofs of \cite[Theorem 1.3]{1} and \cite[Lemma A.2]{nb}, in what follows we prove that similar result is also true to the fractional system (\ref{nb666}) for both $a \neq \beta$ and $a = \beta.$\\
\indent For $(u, v) \in \mathcal{X}_0,$ the energy function of (\ref{nb666}) is given by
\begin{equation}\label{5:limb}
J(u, v)= \frac{1}{2}(\|u\|^{2}_{H^{\frac{1}{2}}} + \|v\|^{2}_{H^{\frac{1}{2}}} )
-\frac{1}{4}\int_{\R^3}\Big(a\phi_{u}u^{2}+ a\phi_{v}v^{2}+2\beta\phi_{u}v^{2}\Big)dx.
\end{equation}
Define
\begin{equation}\label{ciqiong}
c_0 := \inf\Big\{J(u,v): (u, v) \in \mathcal{N}\Big\}, \text{ and }
\end{equation}
$$\mathcal{N} = \Big\{(u, v) \in \mathcal{X}_0 \setminus \{(0, 0)\}: \|u\|^{2}_{H^{\frac{1}{2}}} = \int_{\R^3} (a\phi_{u}u^{2} + \beta\phi_{v}u^{2} )dx, \|v\|^{2}_{H^{\frac{1}{2}}} = \int_{\R^3} (a\phi_{v}v^{2} + \beta\phi_{u}v^{2} )dx,$$
where
$\|u\|^{2}_{H^{\frac{1}{2}}} = \int_{\R^3}\big(|(-\Delta)^{\frac{1}{4}}u|^{2} + u^{2}\big),~\|v\|^{2}_{H^{\frac{1}{2}}} =
\int_{\R^3}\big(|(-\Delta)^{\frac{1}{4}}v|^{2} + v^{2}\big).$
Then, the minimizers of (\ref{ciqiong}) are ground state solutions for (\ref{nb666}).\\
\indent If $(u, v) \in \mathcal{N},$ it is easy to see that
$$J(u, v) = \frac{1}{4}\|u\|^{2}_{H^{\frac{1}{2}}} + \frac{1}{4}\|v\|^{2}_{H^{\frac{1}{2}}}.$$
We claim that $\big(u_0, v_0\big) := \Big(\sqrt{k}Q, \sqrt{l}Q\Big) \in \mathcal{N}$ if $k > 0$ and $l > 0$ satisfying
\begin{equation}\label{4400}
ak + \beta l = 1, \ \text{and} \ al + \beta k =1.\end{equation}
In fact, by the properties (\ref{1:id}) for $Q$, we have
$$\|u_0\|^{2}_{H^{\frac{1}{2}}} = \|\sqrt{k}Q\|^{2}_{H^{\frac{1}{2}}} = 2k\|Q\|^{2} = 2ka^{\ast}, \text{ and } $$
$$\aligned \int_{\R^3}\big(a\phi_{u_0}u_{0}^{2} + \beta\phi_{v_0}u_{0}^{2}\big) &= (ak^2 + \beta kl)\iint_{\R ^3\times\R ^3} \frac{Q^2(x)Q^2(y)}{|x - y|}dxdy = 2k(ak + \beta l)a^{\ast}\\
&= 2ka^{\ast} = \|u_0\|^{2}_{H^{\frac{1}{2}}},~ ~\text{if} ~ak + \beta l = 1.\endaligned$$
Similarly, we know that $\|v_0\|^{2}_{H^{\frac{1}{2}}} = \int_{\R^3}\big(a\phi_{v_0}v_{0}^{2} + \beta\phi_{u_0}v_{0}^{2}\big)dx = 2la^{\ast}$ if $a l+ \beta k = 1$, and then $
\big(u_0, v_0\big) \in \mathcal{N}.$ Therefore, by the definition of (\ref{ciqiong}) we have
\begin{equation}\label{4545}
c_0 \leq J(u_0, v_0) = \frac{1}{4}\|u_0\|^{2}_{H^{\frac{1}{2}}} + \frac{1}{4}\|v_0\|^{2}_{H^{\frac{1}{2}}} = \frac{a^\ast}{2}(k + l).\end{equation}

\begin{lem}\label{lemAf}
If $k > 0, l > 0$ and (\ref{4400}) is satisfied, then $ (u_0, v_0 ) = \big(\sqrt{k}Q, \sqrt{l}Q\big)$ is a ground state solution for \eqref{nb666} provided $a \neq \beta.$
\end{lem}

\noindent{\textbf{Proof.}} By (\ref{4545}), in order to prove this lemma, we need only to show that $J(u_0, v_0) \leq c_0,$ that is, $J(u_0, v_0) = c_0.$ Motivated by \cite{1}, we let $ \{ (u_{n}, v_{n} ) \} \subset \mathcal{X}_0$ be a minimizing sequence for (\ref{ciqiong}), that is,
\begin{equation}\label{4646}
\{ (u_{n}, v_{n} ) \} \subset \mathcal{N} ~\text{and}~ J(u_n, v_n) = \frac{1}{4}\|u_n\|^{2}_{H^{\frac{1}{2}}} + \frac{1}{4}\|v_n\|^{2}_{H^{\frac{1}{2}}} \xrightarrow{n} c_0.
\end{equation}
By Gagliardo-Nirenberg inequality (\ref{GNineq}), we know that
$$\int_{\R^3}\phi_{u_{n}}u_{n}^{2} dx \leq \frac{2}{a^\ast}\|(-\Delta)^{\frac{1}{4}}u_n\|_2^2\|u_n\|_2^2 \leq \frac{1}{2a^\ast}\|u_n\|^{4}_{H^{\frac{1}{2}}}, ~\text{for}~u_n \in H^{\frac{1}{2}}(\R ^3).$$
Then, it follows from $\{ (u_{n}, v_{n} ) \} \subset \mathcal{N}$ that
\begin{equation}\label{4700}
\alpha_{u_n} := \frac{1}{\sqrt{2a^{\ast}}}\Big(\int_{\R^3}\phi_{u_{n}}u_{n}^{2} dx \Big)^\frac{1}{2}  \leq \frac{1}{2a^\ast}\|u_n\|^{2}_{H^{\frac{1}{2}}} = \frac{1}{2a^\ast}\int_{\R^3}\big(a\phi_{u_{n}}u_{n}^{2} + \beta\phi_{v_{n}}u_{n}^{2}\big)dx, \end{equation}
\begin{equation}\label{4800}
\alpha_{v_n} := \frac{1}{\sqrt{2a^{\ast}}}\Big(\int_{\R^3}\phi_{v_{n}}v_{n}^{2} dx\Big)^\frac{1}{2} \leq \frac{1}{2a^\ast}\|v_n\|^{2}_{H^{\frac{1}{2}}} = \frac{1}{2a^\ast}\int_{\R^3}\big(a\phi_{v_{n}}v_{n}^{2} + \beta\phi_{u_{n}}v_{n}^{2}\big)dx. \end{equation}
Moreover, we know from Lemma 2.3 that
\begin{equation}\label{hzy3}\int_{\R^3}\phi_{v_{n}}u_{n}^{2} dx = \int_{\R^3}\phi_{u_{n}}v_{n}^{2} dx =\iint_{\R ^3\times\R ^3} \frac{u^2_{n}(x)v^2_{n}(y)}{|x - y|}dxdy \leq 2a^\ast \alpha_{u_n}\alpha_{v_n}.\end{equation}
Combining (\ref{4700})-(\ref{hzy3}), we have
\begin{equation}\label{hzy5}\begin{split}
&2a^{\ast}\alpha_{u_n} \leq \|u_{n}\|^{2}_{H^{\frac{1}{2}}} \leq 2a^\ast\Big(a\alpha^2_{u_n} + \beta \alpha_{u_n}\alpha_{v_n}\Big),\\
&2a^{\ast}\alpha_{v_n} \leq \|v_{n}\|^{2}_{H^{\frac{1}{2}}} \leq 2a^\ast\Big(a\alpha^2_{v_n} + \beta \alpha_{u_n}\alpha_{v_n}\Big).\end{split}\end{equation}
these and (\ref{4545})-(\ref{4646}) imply that
\begin{equation}\label{411411}
2a^{\ast}\big(\alpha_{u_n} + \alpha_{v_n}\big) \leq \|u_{n}\|^{2}_{H^{\frac{1}{2}}} + \|v_{n}\|^{2}_{H^{\frac{1}{2}}} = 4c_0 + o(1) \leq 2a^{\ast}(k+l) + o(1).\end{equation}
Then, (\ref{hzy5})-(\ref{411411}) give that
\begin{equation}\label{hzy6}
\begin{cases}
\alpha_{u_n} + \alpha_{v_n} \leq k + l + o(1),\\
1 \leq a \alpha_{u_n} + \beta \alpha_{v_n},\\
1 \leq \beta \alpha_{u_n} + a \alpha_{v_n} .\end{cases}\end{equation}
By (\ref{4400}), $ak + \beta l = 1$ and $al + \beta k = 1,$ then, by setting $\gamma_{u_n} = \alpha_{u_n} - k \text{ and } \gamma_{v_n} = \alpha_{v_n} - l $,  (\ref{hzy6}) can be rewritten as
\begin{equation}
\begin{cases}
\gamma_{u_n} + \gamma_{v_n} \leq o(1) \\
a \gamma_{u_n} + \beta \gamma_{v_n} \geq 0,\\
\beta \gamma_{u_n} + a \gamma_{v_n} \geq 0.\end{cases}\end{equation}
Let
$$\Gamma_{n} = \Big\{\big(\gamma_{u_n}, \gamma_{v_n}\big): \gamma_{u_n} + \gamma_{v_n} \leq o(1),a \gamma_{u_n} + \beta \gamma_{v_n} \geq 0,\beta \gamma_{u_n} + a \gamma_{v_n} \geq 0\Big\}.$$
Then, for each $n \geq 1,$ $\Gamma_{n}$ is a triangular region with the vertex $(0, 0)$ and its diameter goes to 0 as $n ~ \rightarrow ~ +\infty$ if $a \neq \beta,$ which then gives that
$$\gamma_{u_n} = \alpha_{u_n} - k \stackrel{n}\rightarrow 0, \gamma_{v_n} = \alpha_{v_n} - l \stackrel{n}\rightarrow 0, ~i.e., \alpha_{u_n} \xrightarrow{n} k, \alpha_{v_n} \xrightarrow{n} l.$$
This and (\ref{411411}) imply that, as $n ~ \rightarrow ~ +\infty,$
$$4c_0 = 2a^{\ast}(k+l), ~i.e.,~c_0 = \frac{a^\ast}{2}(k + l).$$
Hence, $J(u_0, v_0) = \frac{a^\ast}{2}(k + l) =c_0$ by (\ref{4545}), which means that $\big(u_0, v_0\big) = \Big(\sqrt{k}Q, \sqrt{l}Q\Big)$ is a ground state solution. \hfill$\square$
\vskip8pt

\indent In the case of $a = \beta,$ we have the following lemma.
\begin{lem}\label{lemAf2} For any $\theta \in (0, 2\pi),$ $\big(u_0, v_0\big) = \Big(\frac{1}{\sqrt{a}}Q\sin \theta, \frac{1}{\sqrt{a}}Q\cos \theta\Big)$ is a ground state solution for (\ref{nb666}) with $a = \beta.$
\end{lem}

\noindent{\textbf{Proof.}} Motivated by the proof  of Lemma A.2 of \cite{nb},  we consider  the minimization problem
\begin{equation}\label{416}
h:= \inf_{(0,0)\neq (u_{1},u_{2})\in H^{\frac{1}{2}}(\R ^3)\times H^{\frac{1}{2}}(\R ^3)} H\big(u_{1}, u_{2}\big),
\end{equation}
where 
\begin{equation}
H\big(u_{1}, u_{2}\big):= \frac{\int_{\R^3}\Big(|(-\Delta)^{\frac{1}{4}}u_{1}|^{2}+|(-\Delta)^{\frac{1}{4}}u_{2}|^{2}\Big)dx\int_{\R^3}\Big(|u_{1}|^{2}+|u_{2}|^{2}\Big)dx}{a\int_{\R^3}\Big(\phi_{u_1}u_1^2 +\phi_{u_2}u_2^2 + 2\phi_{u_{1}}u_{2}^2\Big)dx}.\end{equation}
By taking $\Big(\frac{1}{\sqrt{2}}Q, \frac{1}{\sqrt{2}}Q\Big)$ as the test function, we have
\begin{equation}\label{zzssll}
h \leq H\Big(\frac{1}{\sqrt{2a}}Q, \frac{1}{\sqrt{2a}}Q\Big) = \frac{\|Q\|_{2}^2}{2a} = \frac{a^{\ast}}{2a}.\end{equation}
On the other hand, by  \cite[Theorem 7.13]{14} and the Gagliardo-Nirenberg inequality (\ref{GNineq}), we have,  for any $(u_{1}, u_{2}) \in H^{\frac{1}{2}}(\R ^3)\times H^{\frac{1}{2}}(\R ^3),$
\begin{equation}\label{zzzz}\begin{split}
H\big(u_{1}, u_{2}\big)&= \frac{\int_{\R^3}\Big(|(-\Delta)^{\frac{1}{4}}u_{1}|^{2}+|(-\Delta)^{\frac{1}{4}}u_{2}|^{2}\Big)dx\int_{\R^3}\Big(|u_{1}|^{2}+|u_{2}|^{2}\Big)dx}{a\int_{\R^3}\Big(\phi_{u_1}u_1^2 +\phi_{u_2}u_2^2 + 2\phi_{u_{1}}u_{2}^2\Big)dx}\\
&\geq \frac{\int_{\R^3}\Big|(-\Delta)^{\frac{1}{4}}\sqrt{|u_{1}|^2 +|u_2|^2}\Big|^{2}dx\int_{\R^3}\Big(\sqrt{|u_{1}|^{2}+|u_{2}|^{2}}\Big)^{2}dx}{a \iint_{\R ^3\times\R ^3} \frac{\Big(\sqrt{|u_{1}|^{2}+|u_{2}|^{2}}\Big)^{2}(x)\Big(\sqrt{|u_{1}|^{2}+|u_{2}|^{2}}\Big)^{2}(y)}{|x - y|}dxdy} \geq \frac{\|Q\|_{2}^2}{2a}.\end{split}\end{equation}
Combining (\ref{zzssll}) and (\ref{zzzz}), we have $h = \frac{\|Q\|_{2}^2}{2a} = \frac{a^{\ast}}{2a}$ and for any $\theta \in (0, 2\pi),$
\begin{equation}\label{422}
\Big(\frac{1}{\sqrt{a}}Q(x)\sin \theta, \frac{1}{\sqrt{a}}Q(x)\cos \theta\Big)\end{equation}
is an optimizer of (\ref{416}). Moreover, one can check that if $\big(u_1, u_2\big)$ is an optimizer of (\ref{416}), then up to scalings, $\big(u_1, u_2\big)$ is a ground state of (\ref{nb666}). Therefore, $\big(u_1, u_2\big)$ being the form of (\ref{422}) is a ground state of
(\ref{nb666}). This completes the proof of this lemma.~~~\hfill$\square$\\

\noindent {\bf Acknowledgements:} This work was supported by  NFSC under Grants No. 11931012, 11871387, 12171379.

\end{document}